\documentclass[10pt]{amsart}

\usepackage{amssymb}
\usepackage{url}
\usepackage{euscript}

\newtheorem{theorem}{Theorem}[section]

\newtheorem{definition}{Definition}[section]
\newtheorem{lemma}{\em Lemma}[section]

\newtheorem{corollary}{\em Corollary}[section]
\newtheorem{remark}{\em Remark}[section]
\newtheorem{question}{\em Question}
\begin{document}

\title[On the preservation of Gibbsianness under symbol amalgamation]
{On the preservation of Gibbsianness\\under symbol amalgamation}
 
\author{Jean--Ren\'e Chazottes \& Edgardo Ugalde}

\begin{abstract}
Starting from the full--shift on a finite alphabet $A$, mingling some symbols 
of $A$, we obtain a new full shift on a smaller alphabet $B$. 
This amalgamation defines a factor map from $(A^{\mathbb N},T_A)$ to $(B^{\mathbb N},T_B)$, 
where $T_A$ and $T_B$ are the respective shift maps.

\noindent
According to the thermodynamic formalism, to each regular function (`potential')
$\psi:A^{\mathbb N}\to{\mathbb R}$, we can associate a unique Gibbs measure $\mu_\psi$. In this article, 
we prove that, for a large class of potentials, the pushforward measure
$\mu_\psi\circ\pi^{-1}$ is still Gibbsian for a potential $\phi:B^{\mathbb N}\to{\mathbb R}$ having 
a `bit less' regularity than $\psi$. In the special case where $\psi$ is a `$2$--symbol' 
potential, the Gibbs measure $\mu_\psi$ is nothing but a Markov measure and the 
amalgamation  $\pi$ defines a hidden Markov chain. In this particular case, our theorem 
can be recast by saying that a hidden Markov chain is a Gibbs measure (for a H\"older potential).
\end{abstract}


\noindent \address{{\em J.-R. Chazottes:} Centre de Physique
   Th\'eorique, CNRS-\'Ecole Polytechnique, 91128 Palaisau Cedex,
   France.}  \email{jeanrene@cpht.polytechnique.fr}

\noindent \address{{\em E. Ugalde:} Instituto de F\'\i sica,
   Univesidad Autonoma de San Luis Potos\'\i, San Luis de Potos\'\i,
   S.L.P., 78290 M\' exico.}  \email{ugalde@ifisica.uaslp.mx}
\date{\today}

\thanks{We acknowledge CONACyT--Mexico, SEP--Mexico, and CNRS--France financial support.}

\maketitle

\tableofcontents

\section{Introduction}~\label{section-introduction}\

\medskip\noindent From different viewpoints and under different names, 
the so--called {\em hidden Markov measures} have received a lot of attention 
in the last fifty years~\cite{karl}.
One considers a (stationary) Markov chain $(X_n)_{n\in{\mathbb N}}$ with finite 
state space $A$ and looks at its `instantaneous' image $Y_n:=\pi(X_n)$, where 
the map $\pi$ is an amalgamation of the elements of $A$ yielding a smaller state 
space, say $B$. It is well--known that in general the resulting chain, 
$(Y_n)_{n\in{\mathbb N}}$, has infinite memory. For concrete examples, see 
{\em e.~g.}~\cite{black} or the more easily accessible reference~\cite{karl} where 
they are recalled.

\medskip\noindent A stationary Markov chain with finite state space $A$ can be 
equivalently defined as a shift--invariant Markov measure $\mu$ on the path space 
$A^{\mathbb N}$ (of infinite sequences of `symbols' from the finite `alphabet' $A$), 
where the shift map $T:A^{\mathbb N}\to A^{\mathbb N}$ is defined by 
$(T\boldsymbol{a})_i=\boldsymbol{a}_{i+1}$. A hidden Markov measure can be 
therefore seen as the the pushforward measure $\mu_\psi\circ\pi^{-1}$ on the path 
space $B^{\mathbb N}$ formed by the instantaneous image under the amalgamation $\pi$, 
of paths in $A^{\mathbb N}$.

\medskip\noindent In the present article, instead of focusing on shift--invariant 
Markov measures, we consider a natural generalization of them. 
Let $\psi:A^{\mathbb N}\to{\mathbb R}$ be a `potential', then, under appropriate 
regularity condition on $\psi$ (see more details below), there is a unique so--called 
Gibbs measure $\mu_\psi$ associated to it. It is a shift--invariant probability measure 
on $A^{\mathbb N}$ with remarkably nice properties. 
Each $r$--step Markov measure falls in this category, since an $r$--step Markov 
measure is nothing but a Gibbs measure defined by a $(r+1)$--symbol potential, 
{\em i.e.}, a potential $\psi$ such that $\psi(\boldsymbol{a})=\psi(\tilde{\boldsymbol{a}})$ 
whenever $\boldsymbol{a}_i=\tilde{\boldsymbol{a}}_i$, $i=0,\ldots,r$, 
with $r$ a strictly positive integer.~\footnote{The case $r=0$ corresponds to product 
measures (i.i.d. process).}
On the other hand, given $\psi$ one can construct a sequence $(\psi_r)$ of 
$(r+1)$--symbol potentials (uniformly approximating $\psi$) such that the sequence of 
associated $r$--step Markov measures $\mu_{\psi_r}$ converges to $\mu_\psi$ (in the 
vague or weak$^*$ topology, at least).

\medskip\noindent Now let $B$ be the alphabet obtained from $A$ by amalgamation of 
some of the symbols of $A$.~\footnote{We assume $B$ has cardinality at least equal 
to two.} The amalgamation defines a surjective ({\em i.e.}, onto) map $\pi:A\to B$ 
which extends to $A^{\mathbb N}$ in the obvious way. 
Given a Gibbs measure $\mu_{\psi}$ on $A^{\mathbb N}$, this map induces a measure 
$\mu_{\psi}\circ\pi^{-1}$ supported on the full shift $B^{\mathbb N}$. 
The question we address now reads:

\begin{question}\label{fullshifts}
Under which condition is the measure $\mu_{\psi}\circ\pi^{-1}$, supported on the 
full shift $B^{\mathbb N}$, still Gibbsian? In other words, under which conditions 
on $\psi$ can one build a `nice' potential $\phi:B^{\mathbb N}\to{\mathbb R}$ such that 
$\mu_{\psi}\circ\pi^{-1}=\mu_\phi$?
In particular, for $\psi$ a $2$--symbol potential, what is the nature of 
$\mu_{\psi}\circ\pi^{-1}$?
\end{question}

\medskip\noindent In this article we make the following answer (made precise below, 
see Theorems~\ref{theorem-projection-of-gibbsian} 
and~\ref{theorem-projection-of-gibbsian-bis}):

\begin{quote}
Under mild regularity condition on $\psi$, the pushforward of the Gibbs measure 
$\mu_\psi$, namely $\mu_{\psi}\circ\pi^{-1}$, is Gibbsian as well, and the associated 
potential $\phi$ can be computed from $\psi$. Furthermore, when $\psi$ is a $2$--symbol 
potential, the corresponding hidden Markov chain is Gibbsian, and it is associated to 
a H\"older potential.
\end{quote}

\medskip\noindent A slightly more general problem is the following. Suppose that we do 
not start with the full shift $A^{\mathbb N}$ but with a {\em subshift of finite type} 
(henceforth SFT) or a {\em topological Markov chain} $X$ \cite{kitchensbook}. The image 
of $X$ is not in general of finite type but it is a {\em sofic subshift}~\cite{kitchensbook}:
\begin{question}\label{sofic}
When $X\subset A^{\mathbb N}$ is a SFT, is the measure $\mu_{\psi}\circ\pi^{-1}$ 
still Gibbsian? 
\end{question}

\medskip\noindent Question~\ref{sofic} has only received very partial answers 
up--to--date. We shall comment on that in Section~\ref{section-concluding-remarks}.

\medskip\noindent The present work is motivated, on the one hand, by our previous 
work in~\cite{2003ChU} in which we attempted to solve Question~\ref{sofic} and were 
partially successful. On the other hand, it was motivated by~\cite{2005ChRU} where 
we were interested in approximating Gibbs measures on sofic subshifts by Markov 
measures on subshifts of finite type. 
Here we combine ideas and techniques both from~\cite{2003ChU} and~\cite{2005ChRU} but 
we need extra work to get more uniformity than previously obtained. 

\medskip\noindent Let us mention two recent works related to ours. 
In \cite{frank}, another kind of transformation of the alphabet is considered,
and the method employed to prove Gibbsianity is completely different from ours.
In \cite{pierre}, the authors study random functions of Markov chains and obtain
results about their loss of memory.

\medskip\noindent The paper is organized as follows.
In the next section we give some notations and definitions. 
In particular, we present the weak${}^*$ convergence of measures as a projective 
convergence and we define the notion of Markov approximants of a Gibbs measure.
In Section~\ref{section-main-result}, we state Theorem~\ref{theorem-projection-of-gibbsian}
which answers Question~\ref{fullshifts} when the starting potential $\psi$ is H\"older 
continuous (its modulus of continuity decays exponentially to $0$).
The proof relies on two lemmas which are proved in Appendices~\ref{proof-of-projections} 
and~\ref{proof-of-projective}, respectively.
In Section~\ref{subexcases}, we generalize Theorem~\ref{theorem-projection-of-gibbsian}
to a class of potentials with subexponential (strictly subexponential or polynomial) decay 
of modulus of continuity.  
We finish (Section~\ref{section-concluding-remarks}) by discussing Question~\ref{sofic} and 
giving a conjecture. Appendix~\ref{refined-perron-frobenius} is devoted to Birkhoff's 
version of Perron--Frobenius theorem for positive matrices, our main tool.

\medskip\noindent We have greatly benefited from the careful reading of an anonymous 
referee whose valuable comments helped us to improve the paper.

\section{Background Material}\


\bigskip\noindent \subsection{Symbolic dynamics}\

\medskip\noindent Let $A$ be a finite set (`alphabet') and $A^{{\mathbb N}}$ be the 
set of infinite sequences of symbols drawn from $A$. We define ${\mathbb N}$ to be 
the set $\{0,1,2,\ldots\}$, that is, the set of positive integers plus $0$.
We denote by $\boldsymbol{a},\boldsymbol{b}$, etc, elements of $A^{\mathbb N}$ and use
the notation $\boldsymbol{a}_m^n$ ($m\leq n$, $m,n \in{\mathbb N}$) for the word 
$\boldsymbol{a}_m\boldsymbol{a}_{m+1}\cdots\boldsymbol{a}_{n-1}\boldsymbol{a}_n$
(of length $n-m+1$). We endow $A^{{\mathbb N}}$ with the distance
\[ d\!_{{\scriptscriptstyle A}}(\boldsymbol{a},\boldsymbol{b})
      :=\left\{ \begin{array}{l}    
                \exp \left( -\min\{n\geq 0:\boldsymbol{a}_0^n\neq \boldsymbol{b}_0^n\}
                     \right)         \quad\textup{if}\; \boldsymbol{a}\neq \boldsymbol{b} \\
                     0               \quad\textup{otherwise}. \end{array} \right.
\]
The resulting metric space $(A^{\mathbb N},d\!_{{\scriptscriptstyle A}})$ is compact.

\medskip\noindent The {\em shift transformation} $T:A^{\mathbb N}\to A^{\mathbb N}$ is 
defined by $(T\boldsymbol{a})_n=\boldsymbol{a}_{n+1}$ for all $n\in{\mathbb N}$.

\medskip\noindent A {\em subshift} $X$ of $A^{\mathbb N}$ is a closed $T$--invariant 
subset of $A^{\mathbb N}$. 

\medskip\noindent Given a set of {\em admissible words} ${\mathcal L}\subset A^\ell$ for 
some fixed integer $\ell\geq 2$, one defines a {\em subshift of finite type} 
$A_{{\mathcal L}}\subset A^{\mathbb N}$ by
\[
A_{{\mathcal L}}:=\{\boldsymbol{a}\in A^{\mathbb N}:\ 
\boldsymbol{a}_n^{n+\ell-1}\in{\mathcal L} , \ \forall n\in{\mathbb N}\}.
\]
A subshift of finite type defined by words in ${\mathcal L}\subset A^2$ is called a
{\em topological Markov chain}.
It can be equivalently described by the transition matrix $M:A\times A\to\{0,1\}$
such that $M(a,b)=\chi_{{\scriptstyle{\mathcal L}}}(ab)$, where 
$\chi_{{\scriptstyle{\mathcal L}}}$ is the indicator function of the set ${\mathcal L}$.
We will use both $A_{{\mathcal L}}$ and $A_M$ to denote the corresponding subshift of 
finite type.

\medskip\noindent Note that the `full shift' $(A^{\mathbb N},T)$ can be seen as the 
subshift of finite type defined by all the words of length $\ell$, and we have the 
identification $A^{\mathbb N}\equiv A_{A^\ell}$.

\medskip\noindent Let $X\subset A^{\mathbb N}$ be a subshift.
A point $\boldsymbol{a}\in X$ is {\em periodic} with period $p\geq 1$ if 
$T^p\boldsymbol{a}=\boldsymbol{a}$, and $p$
is its minimal period if in addition $T^k\boldsymbol{a}\neq \boldsymbol{a}$ whenever
$0 < k < p$.
We denote by $\hbox{\rm Per}_p(X)$ the collection of all periodic points with 
period $p$ in $X$, and by $\hbox{\rm Per}(X)$ the collection of all periodic points 
in $X$, {\em i.e.}, $\hbox{\rm Per}(X)=\cup_{p\geq 1}\hbox{\rm Per}_p(X)$.

\medskip\noindent Given an arbitrary subshift $X\subset A^{\mathbb N}$ and $m\in {\mathbb N}$, 
the set of all the {\em $X$--admissible words of length $m+1$} is  the set
\[X_m:=\{\boldsymbol{w}\in A^{m+1}:\ \exists\ \boldsymbol{a}\in X, \ 
           \boldsymbol{w}=\boldsymbol{a}_0^m\}. \]

\medskip\noindent
It is a well known fact that a topological Markov chain $A_M$ is
topologically mixing if its transition matrix $M$ is primitive, {\em  i.e.},
if and only if there exists an integer $n\geq 1$ such that 
$M^{n}>0$~\footnote{On the other hand, if none of the rows or columns of $M$ is identically 
zero, $A_M$ is topologically mixing implies $M$ is primitive.}. In this case, the smallest 
of such integers is the so called primitivity index of $M$.

\medskip\noindent For a subshift $X\subset A^{\mathbb N}$, $\boldsymbol{w}\in X_m$ and 
$m\in{\mathbb N}$, the set
\[ [\boldsymbol{w}]:=\{\boldsymbol{a}\in X:\ \boldsymbol{a}_0^{m}=\boldsymbol{w}\} \]
is the {\em cylinder} based on $\boldsymbol{w}$.

\medskip\noindent We will use boldfaced symbols $\boldsymbol{a},\boldsymbol{b}$, etc, not 
only for infinite sequences but also for finite ones ({\em i.e.}, for words). The context 
will make clear whether we deal with a finite or an infinite sequence.


\bigskip\noindent\subsection{Thermodynamic formalism}\

\medskip\noindent For a subshift $X\subset A^{{\mathbb N}}$, cylinders are clopen sets 
and generate the Borel $\sigma$--algebra. We denote by $\EuScript{M}(X)$ the set of Borel 
probability measures on $X$ and by {\em $\EuScript{M}_T(X)$ the subset of $T$--invariant 
probability measures} on $X$. Both are compact convex sets in weak$^*$ topology.
The weak$^*$ topology can be metrized \cite{bowenbook} by the distance 
\[ D(\mu,\nu):= \sum_{m=0}^\infty 2^{-(m+1)} \left(\sum_{\boldsymbol{w}\in X_m}
  \left| \mu[\boldsymbol{w}]-\nu[\boldsymbol{w}]\right|\right). \]

\medskip\noindent It turns out that the following
notion of convergence is very convenient in our later calculations.

\begin{definition}
We say that a sequence $(\mu_n)_{n\in{\mathbb N}}$ of probability
measures in $\EuScript{M}(X)$ {\em converges in the projective sense} to a
measure $\mu\in \EuScript{M}(X)$ if for all $\epsilon > 0$ and $N >1$ there exists 
$N'>1$ such that
\[ \exp(-\epsilon)\leq\frac{\mu_n[\boldsymbol{w}]}{\mu[\boldsymbol{w}]}\leq\exp(\epsilon) \]
for all admissible words $\boldsymbol{w}$ of length $k\leq N$, and for all $n\geq N'$.
\end{definition} 

\medskip\noindent It is easy to verify that convergence in the projective
sense implies weak${}^*$ convergence. On the other hand, when all the 
measures involved share the same support, weak${}^*$ and projective 
convergence coincide. Though it is the case in this paper, we will speak of
projective convergence. 

\medskip\noindent We make the following definitions.

\begin{definition}[$(r+1)$--symbol potentials]
A function $\psi:A^{\mathbb N}\to{\mathbb R}$ will be called a {\em potential}. 
We say that a potential $\psi:A^{\mathbb N}\to {\mathbb R}$ is an $(r+1)$--symbol potential 
if there is an $r\in{\mathbb N}$ such that  
\[ \psi(\boldsymbol{a})=\psi(\boldsymbol{b}) \text{ whenever } 
                      \boldsymbol{a}_0^r=\boldsymbol{b}_0^r. \]
Of course, we take $r$ to be the smallest integer with this property. 

\noindent 
We will say that $\psi$ is {\em locally constant} if it is an $(r+1)$--symbol potential 
for some $r\in{\mathbb N}$. 
\end{definition}

\medskip\noindent A way of quantifying the regularity of a potential 
$\psi:A^{\mathbb N}\to {\mathbb R}$ is by using its modulus of continuity on cylinders, 
or variation, defined by
\[ \hbox{\rm var}_n\psi:=\sup\{|\psi(\boldsymbol{a})-\psi(\boldsymbol{b})|:
\ \boldsymbol{a},\boldsymbol{b}\in A^{\mathbb N}, \ \boldsymbol{a}_0^n=\boldsymbol{b}_0^n\}.\]

\medskip\noindent A potential $\psi$ is continuous if and only if 
$\hbox{\rm var}_n\psi\to 0$ as $n\to\infty$. An $(r+1)$--symbol potential $\psi$ can 
be alternatively defined by requiring that $\hbox{\rm var}_n\psi=0$ whenever $n\geq r$, 
and thus it is trivially continuous. If there are ${\mathcal C}>0$ and 
$\varrho\in]0,1[$ such that $\hbox{\rm var}_n\phi\leq {\mathcal C}\varrho^n$ 
for all $n\geq 0$, then $\psi$ is said to be H{\"o}lder continuous. 

\medskip\noindent We will use the notation
\[
S_n\psi(\boldsymbol{a}):=\sum_{k=0}^{n-1}\psi\circ T^k(\boldsymbol{a}), \ n=1,2,\ldots
\]
Throughout we will write 
\[
x\lessgtr y\, C^{\pm 1}\quad\textup{for}\quad  C^{-1}\leq \frac{x}{y}\leq C 
\]
for $x,y$ and $C$ strictly positive numbers. Accordingly we will use the
notation $x\lessgtr y\exp(\pm C)$. We also write $x\lessgtr y \pm C$ for
$-C\leq x-y\leq C$.

\medskip\noindent We now define the notion of Gibbs measure we will use in the sequel. 

\begin{definition}[Gibbs measures]
Let $X\subset A^{\mathbb N}$ be a subshift and $\psi:A^{\mathbb N}\to{\mathbb R}$ be a 
potential such that $\psi|_{X}$ is continuous. A measure $\mu\in\EuScript{M}_T(X)$ is 
a {\em Gibbs measure for the potential $\psi$},  if there are constants $C=C(\psi,X) \geq 1$ 
and $P=P(\psi,X)\in {\mathbb R}$ such that
\begin{equation}~\label{gibbs-inequality}
\frac{\mu[\boldsymbol{a}_0^n]}{\exp(S_{n+1}\psi(\boldsymbol{a})-(n+1)P)} 
\lessgtr C^{\pm 1},
\end{equation}
for all $n\in {\mathbb N}$ and $\boldsymbol{a}\in X$. We denote by $\mu_\psi$ such a measure.
\end{definition}

\medskip\noindent The constant $P=P(\psi,X)$ is the {\em topological pressure}~\cite{kellerbook} 
of $X$ with respect to $\psi$. It can be obtained, for $X$ a subshfit of finite type, as 
follows:
\begin{equation}~\label{topological-pressure}
P(\psi,X)=\limsup_{n\to\infty}\frac{1}{n}\log 
\sum_{\boldsymbol{a}\in\hbox{\rm \tiny Per}_{n}(X)}\exp(S_n\psi(\boldsymbol{a})).
\end{equation}
We will say that {\em the potential $\psi$ is normalized on $X$} if $P(\psi,X)=0$. We can 
always normalize a potential $\psi$ by replacing $\psi$ by $\psi-P(\psi,X)$. This does 
not affect the associated Gibbs measure $\mu_\psi$.

\medskip\noindent In the above definition, we allow that $\psi=-\infty$ on 
$A^{\mathbb N}\backslash X$. In other words, $\psi$ is upper semi--continuous on 
$A^{\mathbb N}$. 

\begin{remark} 
If $\mu\in \EuScript{M}_T(X)$ is such that the sequence
$\left(\log (\mu[\boldsymbol{a}_0^n]/\mu[\boldsymbol{a}_1^n])\right)_{n=1}^\infty$ converges 
uniformly in $\boldsymbol{a}\in X$, then the potential $\psi: X\to {\mathbb R}$ given by
\begin{equation}\label{psibymu}
\psi(\boldsymbol{a})=\lim_{n\to\infty} \log 
\left(\frac{\mu[\boldsymbol{a}_0^n]}{\mu[\boldsymbol{a}_1^n]}\right)
\end{equation}
is continuous on $X$, and $\mu$ is a Gibbs measure with respect to $\psi$, {\em i.e.}
$\mu=\mu_\psi$. Furthermore, $\psi$ is such that $P(\psi)=0$.

\medskip\noindent Notice that $\mu[\boldsymbol{a}_0^n]/\mu[\boldsymbol{a}_1^n]$
is nothing but the probability under $\mu$ of $\boldsymbol{a}_0$ given $\boldsymbol{a}_1^{n}$. Therefore,
by the martingale convergence theorem  
the sequence 
$\left(\log (\mu[\boldsymbol{a}_0^n]/\mu[\boldsymbol{a}_1^n])\right)_{n=1}^\infty$ converges 
for $\mu$--a.~e. $\boldsymbol{a}\in X$. The uniform convergence is what makes $\mu$ a Gibbs 
measure.
\end{remark}

\medskip\noindent We have the following classical theorem.

\begin{theorem}[\cite{lr}]\label{scu}
Let $X\subset A^{\mathbb N}$ be a topologically mixing subshift of finite type and 
$\psi:X\to{\mathbb R}$. If 
\begin{equation}\label{sumvar}
\sum_{n=0}^\infty\hbox{\rm var}_n\psi<\infty
\end{equation}
then there exists a unique Gibbs measure $\mu_\psi$, {\em i.e.}, a unique 
$T$-invariant probability measure satisfying~\eqref{gibbs-inequality}.
\end{theorem}

\medskip\noindent
\begin{remark}
By this theorem we have a partial converse to~\eqref{psibymu} in the sense that
there the potential is defined by the measure, while in the theorem it is the 
potential which defines the measure.

\medskip\noindent Notice that the uniqueness part of the theorem is granted by the Gibbs 
inequality~\eqref{gibbs-inequality}, since two measures satisfying it have to be 
absolutely continuous with respect to each other. It is the existence part which is
nontrivial.
\end{remark}

\medskip\noindent For a proof of Theorem~\ref{scu} see {\em e.g.}~\cite{kellerbook}. This 
includes the case of H\"older continuous potentials treated in, {\em e.g.}, 
\cite{bowenbook,ruellebook}.

\bigskip\noindent\subsection{Markov measures and Markov approximants} \

\medskip\noindent Markov measures can be seen as Gibbs measures.
Colloquially, an $r$--step Markov measure is defined by the property
that the probability that $\boldsymbol{a}_n=a\in A$ given $\boldsymbol{a}_0^{n-1}$ 
depends only on $\boldsymbol{a}_{n-r}^{n-1}$.~\footnote{We assume that $r\geq 1$. 
The case $r=0$ corresponds to an i.i.d. process, in which case the Gibbsianity is evident.} 
What is usually called a Markov measure corresponds to $1$--step 
Markov measures. On the full shift, the case $r=0$ gives product measures.
A $T$--invariant probability measure is an $r$--step Markov if and only if it is the 
Gibbs measure of an $(r+1)$--symbol potential.
Given an $(r+1)$--symbol potential $\psi$, which we identify as a function on $A^{r+1}$,
one can define the transition matrix
${\mathcal M}_\psi: A^r\times A^r\to{\mathbb R}^+$ such that
\[
{\mathcal M}_{\psi}(\boldsymbol{v},\boldsymbol{w}):=\left\{\begin{array}{ll}
\exp(\psi(\boldsymbol{v}\boldsymbol{w}_{r-1}))
& \textup{ if  } \boldsymbol{v}_1^{r-1}=\boldsymbol{w}_0^{r-2}, \\
0 & \text{ otherwise.}\end{array}\right.
\]
By $\boldsymbol{v}\boldsymbol{w}_{r-1}$ we mean the word obtained by
concatenation of $\boldsymbol{v}$ and $\boldsymbol{w}_{r-1}$ (the last letter of 
$\boldsymbol{w}$).

\medskip\noindent By Perron--Frobenius Theorem (cf. 
Appendix~\ref{refined-perron-frobenius}) there exist a right eigenvector 
$\bar{R}_\psi>0$ such that $\sum_{\boldsymbol{a}\in A^r}\bar{R}_\psi(\boldsymbol{a})=1$, 
and a left eigenvector $\bar{L}_\psi>0$ such that $\bar{L}_\psi^{\dag}\bar{R}_\psi=1$,
associated to the maximal eigenvalue $0<\rho_\psi:=\max{\rm spec}({\mathcal M}_\psi)$.
Then the measure $\mu$ defined by
\begin{equation}\label{parry-measure}
\mu[\boldsymbol{a}_0^n]:=\bar{L}_\psi(\boldsymbol{a}_0^{r-1})\ 
\frac{
\prod_{j=0}^{n-r}{\mathcal M}_\psi(\boldsymbol{a}_j^{j+r-1},\boldsymbol{a}_{j+1}^{j+r})
    }{\rho_\psi^{n-r+1}}\ \bar{R}_\psi(\boldsymbol{a}_{n-r+1}^{n}),
\end{equation}
for each $\boldsymbol{a}\in A^{\mathbb N}$ and $n\in {\mathbb N}$ such that $n\geq r$, 
is easily seen to be a $T$--invariant probability measure satisfying~\eqref{gibbs-inequality} 
with
\[ P=\log(\rho_\psi)\quad\textup{and}\quad
   C=\rho_\psi^r\,e^{-r\|\psi\|}\, 
   \frac{\displaystyle \max\{\bar{L}_\psi(\boldsymbol{w})\bar{R}_\psi(\boldsymbol{w}'):\ 
    \boldsymbol{w},\boldsymbol{w}'\in A^r\}
       }{\displaystyle \min\{\bar{L}_\psi(\boldsymbol{w})\bar{R}_\psi(\boldsymbol{w}'):\ 
       \boldsymbol{w},\boldsymbol{w}'\in A^r\}},
\]
where $\|\psi\|:=\sup\{|\psi(\boldsymbol{a})|:\ \boldsymbol{a}\in A^{\mathbb N}\}$.
Therefore $\mu=\mu_\psi$ is the unique Gibbs measure associated to the $(r+1)$--symbol 
potential $\psi$.

\medskip\noindent {\bf Markov and locally constant approximants}. \

\medskip\noindent Given a continuous $\psi:A^{\mathbb N}\to{\mathbb R}$, one can uniformly approximate it 
by a sequence of $(r+1)$--symbol potentials $\psi_r$, $r=1, 2,\ldots$, in such a way that 
$\|\psi-\psi_r\|\leq \hbox{\rm var}_r(\psi)$, which goes to $0$ as $r$ goes to $\infty$.
The $\psi_r$'s are not defined in a unique way but this does not matter since 
the associated $r$--step Markov measures $\mu_{\psi_r}$, approximate the same 
Gibbs measure $\mu_\psi$.
We can choose $\psi_r(\boldsymbol{a}):=\max\{\psi(\boldsymbol{b}):\ 
\boldsymbol{b}\in[\boldsymbol{a}_0^r]\}$ for instance. 

\medskip\noindent 
The potential $\psi_r$ will be called the {\em (r+1)--symbol approximant} of 
$\psi$ and the associated $r$--step Markov measure $\mu_{\psi_r}$ will be
the {\em rth Markov approximant} of $\mu_\psi$. It is well known 
(and not difficult to prove) that $\mu_{\psi_r}$ converges in the weak$^*$ 
topology to $\mu_\psi$.

\medskip\noindent \section{Main Result}~\label{section-main-result}\

\medskip\noindent 
The next theorem answers Question \ref{fullshifts} when $\psi$ is H\"older continuous
(Theorem~\ref{theorem-projection-of-gibbsian}).
For the sake of simplicity we discuss the generalization of that theorem to
a class of less regular potentials ({\em i.e.}, $\hbox{\rm var}_n(\psi)$ decreases 
subexponentially or polynomially) in Section \ref{subexcases}.

\medskip\noindent {\bf Amalgamation map}. 
Let $A, B$ be two finite alphabets, with $\textup{Card}(A) >\textup{Card}(B)$,
and $\pi:A\to B$ be a surjective map (`amalgamation') which extends 
to the map $\pi:A^{\mathbb N}\to B^{\mathbb N}$ (we use the same letter for both) such that 
$(\pi\boldsymbol{a})_n=\pi(\boldsymbol{a}_n)$ for all $n\in {\mathbb N}$. The map $\pi$ 
is continuous and shift--commuting, {\em i.e.}, it is a factor map from $A^{\mathbb N}$ 
onto $B^{\mathbb N}$.

\medskip\noindent
\begin{theorem}\label{theorem-projection-of-gibbsian}
Let $\pi:A^{\mathbb N}\to B^{\mathbb N}$ be the amalgamation map just defined and 
$\psi:A^{\mathbb N}\to{\mathbb R}$ be a H\"older continuous potential.
Then the measure $\mu_{\psi}\circ\pi^{-1}$ is a Gibbs measure with support $B^{\mathbb N}$,
for a potential $\phi:B^{\mathbb N}\to{\mathbb R}$ such that 
\[
\hbox{\rm var}_n(\phi)\leq {\mathcal D} \exp(-c \sqrt{n})
\]
for some $c,{\mathcal D}>0$, and all $n\in {\mathbb N}$. 

\noindent
Furthermore, this potential $\phi:B^{\mathbb N}\to{\mathbb R}$ is normalized and it is 
given by
\begin{equation}\label{formulepourphi}
\phi(\boldsymbol{b})=\lim_{r\to\infty}\lim_{n\to\infty}
 \log\left(\frac{\mu_{\psi_r}\circ\pi^{-1}[\boldsymbol{b}_0^n]
               }{\mu_{\psi_r}\circ\pi^{-1}[\boldsymbol{b}_1^n]}\right),
\end{equation}
where $\psi_r$ is the $(r+1)$--symbol approximant of $\psi$.

\medskip\noindent If $\psi$ is locally constant, then for all $n$
\[
\hbox{\rm var}_n(\phi)\leq C \vartheta^{n}
\] 
where $\vartheta\in ]0,1[$, $C>0$.
\end{theorem}

\medskip\noindent The case of locally constant potentials in the theorem can be rephrased 
as follows:

\noindent 
When $\mu_\psi$ is an $r$--step Markov measure, with $r > 0$, the pushforward measure
$\mu_{\psi}\circ\pi^{-1}$, {\em i.e.} the hidden Markov measure, is a Gibbs measure for a 
H\"older continuous potential $\phi$ given by 
\begin{equation}\label{phibynu}
\phi(\boldsymbol{b})=\lim_{n\to\infty} 
\log \left(\frac{\mu_\psi\circ\pi^{-1}[\boldsymbol{b}_0^n]
               }{\mu_\psi\circ\pi^{-1}[\boldsymbol{b}_1^n]}\right).
\end{equation}
The case $r=0$ is trivial: the Gibbs measure is simply a product measure and its pushforward
is also a product measure.

\medskip\noindent The proof of Theorem \ref{theorem-projection-of-gibbsian} relies on the 
following two lemmas whose proofs are deferred to Appendices \ref{proof-of-projections} and 
\ref{proof-of-projective}.

\medskip\noindent
\begin{lemma}[Amalgamation for $(r+1)$--symbol potentials]\label{lemma-projections}
The measure $\mu_{\psi_r}\circ\pi^{-1}$, with $r >0$, is a Gibbs measure for the potential 
$\phi_r:B^{\mathbb N}\to{\mathbb R}$ obtained as the following limit
\begin{equation}\label{lephir}
\phi_r(\boldsymbol{b}):=\lim_{n\to\infty}
\log\left(\frac{\mu_{\psi_r}\circ\pi^{-1}[\boldsymbol{b}_0^n]
               }{\mu_{\psi_r}\circ\pi^{-1}[\boldsymbol{b}_1^n]}\right).
\end{equation}
Furthermore, there are constants $C>0$ and $\theta\in [0,1[$ such that, for any
positive integer $n > r$ and for any $\boldsymbol{b}\in B^{\mathbb N}$ we have
\begin{equation}\label{regularitephir}
\left| \phi_r(\boldsymbol{b})-
\log\left(\frac{\mu_{\psi_r}\circ\pi^{-1}[\boldsymbol{b}_0^n]
              }{\mu_{\psi_r}\circ\pi^{-1}[\boldsymbol{b}_1^n]}\right)\right|\leq C\, 
              r^2\, \theta^{\frac{n}{r}}.
\end{equation}
\end{lemma}

\medskip\noindent
\begin{lemma}[Projective convergence of Markov approximants]
\label{lemma-projective}
The sequence of measures $(\mu_{\psi_{r}})$ converges in the projective sense to the Gibbs 
measure $\mu_\psi$ associated to the potential $\psi$. 

\noindent Furthermore, for all $n,\,r >0$ and $\boldsymbol{w}\in A^n$, we have 
$\mu_{\psi_r}[\boldsymbol{w}]\lessgtr \mu_\psi[\boldsymbol{w}] \,\exp(\pm \epsilon_{r,n})$, 
where
\begin{equation}\label{epsilonr}
\epsilon_{r,n}:= D\,\sum_{s=r}^{\infty}\left( (n+ (s+1)(s+2)) 
               \hbox{\rm var}_s\psi+ s\,\theta^{s}\right),
\end{equation}
for adequate constants $D > 0$ and $\theta\in[0,1[$ (the same $\theta$ as in 
Lemma~\ref{lemma-projections}).
\end{lemma}

\medskip\noindent With the two previous lemmas at hand, we can proceed to the proof of 
Theorem~\ref{theorem-projection-of-gibbsian}.

\bigskip\noindent {\em Proof of Theorem \ref{theorem-projection-of-gibbsian}}\

\noindent We start by proving that the sequence $(\mu_{\psi_r}\circ \pi^{-1})_r$ converges 
in the projective sense to $\mu_\psi\circ\pi^{-1}$.

\noindent On the one hand, Lemma~\ref{lemma-projections} tells us that the measure
$\nu_r:=\mu_{\psi_r}\circ\pi^{-1}$
is Gibbsian for the potential $\phi_r:B^{\mathbb N}\to{\mathbb R}$ given by
\[
\phi_r(\boldsymbol{b})=\lim_{n\to\infty}
\log\left(\frac{\nu_r[\boldsymbol{b}_0^n]}{\nu_r[\boldsymbol{b}_1^n]}\right).
\]

\medskip\noindent On the other hand, Lemma~\ref{lemma-projective} ensures that for
each $n,\,r >0$ with $n\geq r$, and each $\boldsymbol{v}\in A^{n}$, we have 
$\mu_{\psi_r}[\boldsymbol{v}]\lessgtr \mu_\psi[\boldsymbol{v}] \exp(\pm\,\epsilon_{r,n})$ 
where $\epsilon_{r,n}$ is defined as in~\eqref{epsilonr}. From this it follows that 
for each $\boldsymbol{w}\in B^{n}$ we have
\begin{eqnarray}~\label{projective-nu}
\nu_r[\boldsymbol{w}]&:=& 
\sum_{\boldsymbol{v}{\scriptscriptstyle \in A^{n}}:\pi\boldsymbol{v}=\boldsymbol{w}} 
     \mu_{\psi_r}[\boldsymbol{v}]\nonumber\\
               &\lessgtr& 
\exp(\pm\epsilon_{r,n})\, 
\sum_{\boldsymbol{v}{\scriptscriptstyle\in A^{n}}:\pi\boldsymbol{v}=\boldsymbol{w}} 
     \mu_{\psi}[\boldsymbol{v}]\nonumber\\
         &\lessgtr& \exp(\pm\epsilon_{r,n})\, \mu_\psi\circ\pi^{-1}[\boldsymbol{w}].
\end{eqnarray}
Otherwise said, the sequence of approximants $(\nu_r\equiv \mu_{\psi_r}\circ\pi^{-1})_{r}$,
converges in the projective sense to the induced measure $\mu_{\psi}\circ\pi^{-1}$, 
and the speed of convergence is the same both the factor and the original system.

\medskip\noindent Now we prove that the pushforward measure $\nu:=\mu_\psi\circ\pi^{-1}$ 
is a Gibbs measure. 

\noindent According to Lemma~\ref{lemma-projections} and Eq.~\eqref{projective-nu}, for 
any $\boldsymbol{b}\in B^{\mathbb N}$, and $n,\, r >0$ with $n \geq r$, we have
\begin{equation}~\label{induced-potential}
\left|\phi_r(\boldsymbol{b})-
\log\left(\frac{\nu[\boldsymbol{b}_0^n]
              }{\nu[\boldsymbol{b}_1^n]}\right)\right|\leq 2\,
           \epsilon_{r,n}+ C\,r^2\,\theta^{\frac{n}{r}}.
\end{equation}
Let us take, for each $r > 0$, $n=n(r):=r^2$, and let $r^* > 0$ be such that both 
$s\mapsto s^2\theta^s$ and $s\mapsto \epsilon_{s,s^2}$ define decreasing functions 
in $[r^*,\infty)$. Hence, using the triangle inequality we obtain 
\[|\phi_r(\boldsymbol{b})-\phi_{r'}(\boldsymbol{b})|
\leq 2\left(2\epsilon_{r,r^2}+C\,r^2\theta^{r}\right)
\] 
for all $r^*\leq r <r'$, and for any $\boldsymbol{b}\in B^{\mathbb N}$. This proves uniform convergence of the 
sequence of potentials $(\phi_r)_{r}$. The limit is the continuous function 
$\phi:B^{\mathbb N}\to{\mathbb R}$ defined by
\[
\phi(\boldsymbol{b}):=\lim_{n\to\infty}
\log\left(\frac{\nu[\boldsymbol{b}_0^n]}{\nu[\boldsymbol{b}_1^n]}\right).
\]

\medskip\noindent If we verify that $\phi$ satisfies condition~\eqref{sumvar}, then, 
according to the 
observation following Theorem~\ref{scu}, this will prove that 
$\nu\equiv \mu_\psi\circ\pi^{-1}$ is the unique Gibbs measure for $\phi$. From 
Ineq.~\eqref{induced-potential} it follows that
\begin{eqnarray*}
|\phi(\boldsymbol{b})-\phi(\tilde{\boldsymbol{b}})| & \leq & 
 \big|\phi(\boldsymbol{b})-\phi_{r}(\boldsymbol{b})\big| + \left|\phi_{r}(\boldsymbol{b})-
 \log\left(\frac{\nu[\boldsymbol{b}_0^n]}{\nu[\boldsymbol{b}_1^n]}\right)\right|\\
                                                    &     &
  +\left|\phi_{r}(\tilde{\boldsymbol{b}})-
  \log\left(\frac{\nu[\boldsymbol{b}_0^n]}{\nu[\boldsymbol{b}_1^n]}\right)\right|+
                   |\phi_{r}(\tilde{\boldsymbol{b}})-\phi(\tilde{\boldsymbol{b}})|\\
                                                   &\leq & 
 4\left(2\epsilon_{r,r^2}+C\,r^2\theta^{r}\right)+
                           2\left(2\epsilon_{r,n}+C\,r^2\,\theta^{\frac{n}{r}}\right),
\end{eqnarray*}
for all $\boldsymbol{b},\tilde{\boldsymbol{b}}\in B^{{\mathbb N}}$ such that 
$\tilde{\boldsymbol{b}}\in [\boldsymbol{b}_0^n]$, and every $n > r\geq r^*$.

\noindent Since $\psi$ is  H\"older continuous and 
\[\epsilon_{r,r^2}:=D\sum_{s=r}^{\infty}(r^2+(s+1)(s+2))\hbox{\rm var}_s\psi+s\theta^s),\] 
then there exist
${\mathcal C} > 0$ and $\varrho\in[\theta,1[$ (remember that $\theta\in [0,1[$) such that 
$\max(\epsilon_{r,r^2},r^2\theta^r)\leq {\mathcal C}\varrho^r$.
We take again $n=n(r)=r^2$ and obtain, for all $n\in{\mathbb N}$, 
\[ \hbox{\rm var}_n\phi\leq {\mathcal D} \exp(-c \sqrt{n}) \]
with a convenient ${\mathcal D}\geq 6\,{\mathcal C}(2+C)$, and $c=-\log(\varrho)$.

\medskip\noindent The case of a locally constant $\psi$ is the immediate consequence of 
Lemma~\ref{lemma-projections} and one has $\vartheta=\theta^{\frac1r}$.

\medskip\noindent The theorem is now proved $\qed$

\medskip\noindent \begin{remark}
The competition between the terms $\epsilon_{r,n}$ and $\theta^{n/r}$ in the upper 
bound of $\hbox{\rm var}_n\phi$ leads to a subexponential bound, namely
$\hbox{\rm var}_n\phi\leq {\mathcal D} \exp(-c n^{\frac{\delta}{1+\delta}})$, for any 
$\delta>0$. We made the choice $\delta=1$.
\end{remark}

\medskip\noindent \section{Generalization to less regular potentials}\label{subexcases}

\medskip\noindent In this section we go beyond H\"older continuous potentials and look at 
potentials $\psi$ such that $\hbox{\rm var}_r(\psi)$ decreases slower than exponentially.
Besides the fact that $\sum_r\hbox{\rm var}_r\psi < \infty$ is always assumed, the only place
where a finer control in the  decrease of $\hbox{\rm var}_r(\psi)$ is required, is 
inside the proof of Lemma~\ref{lemma-projective}. There, the projective convergence of the 
Markov approximants depends on the fact that
\[ \epsilon_{r,n}:=D\sum_{s=r}^\infty ((n+ (s+1)(s+2)) 
        \hbox{\rm var}_s\psi+s\theta^s)\rightarrow 0,
       \text{ when } r\rightarrow \infty,\]
for each $n > 0$. Furthermore, the variation of the induced potential, $\hbox{\rm var}_n\phi$, 
is upper bounded by a linear combination of $\epsilon_{\sqrt{n},n}$ and $n\, \theta^{n/r}$. 
After this consideration, we can 
generalize Theorem~\ref{theorem-projection-of-gibbsian} as follows.

\medskip\noindent
\begin{theorem}\label{theorem-projection-of-gibbsian-bis}
Let $\pi:A^{\mathbb N}\to B^{\mathbb N}$ be the amalgamation map just defined and 
$\psi:A^{\mathbb N}\to{\mathbb R}$ be such that 
$\sum_{s=0}^\infty s^2\,\hbox{\rm var}_s\psi<\infty$.
Then the measure $\mu_{\psi}\circ\pi^{-1}$ is a Gibbs measure with support $B^{\mathbb N}$
for a normalized potential $\phi:B^{\mathbb N}\to{\mathbb R}$ defined by the limit
\[\phi(\boldsymbol{b})=\lim_{r\to\infty}\lim_{n\to\infty}
 \log\left(\frac{\mu_{\psi_r}\circ\pi^{-1}[\boldsymbol{b}_0^n]
               }{\mu_{\psi_r}\circ\pi^{-1}[\boldsymbol{b}_1^n]}\right), \]
where $\psi_r$ is the $(r+1)$--symbol approximant of $\psi$.

\medskip\noindent \underline{If $\hbox{\rm var}_n\psi$ has subexponential decreasing}, 
{\em i.e.}, if $\hbox{\rm var}_n\psi \leq C \exp(-c\,n^{\gamma})$ for some $c,\,C >0$ 
and $\gamma\in]0,1[$, then there are constants $D > C$ and $0 < d < c$ such that
\[
\hbox{\rm var}_n(\phi)\leq D\exp\left(-d\, n^{\frac{\gamma}{1+\gamma}}\right)\]
for all $n\in{\mathbb N}$. 

\medskip\noindent \underline{If $\hbox{\rm var}_n\psi$ is polynomially decreasing}, 
{\em i.e.}, if $\hbox{\rm var}_n\psi\leq C n^{-q}$, for some $C >0$ and $q > 3$, then 
for all $\epsilon \in (0,q-3)$ there is a constant $D > C$ such that
\[\hbox{\rm var}_n(\phi)\leq D\,\frac{1}{n^{q-2-\epsilon}}\]
for all $n\in{\mathbb N}$. 
\end{theorem}

\medskip\noindent 
\begin{remark} As mentioned above, the $n$--variation of the induced potential is upper 
bounded by linear combination of $\epsilon_{r,n}$ and $r^2\theta^{n/r}$. 
We have to optimize the choice of the function $r\mapsto n(r)$ in such a way that 
$n/r\rightarrow \infty$ when $r\to\infty$, and that the resulting $n$--variation of $\psi$ 
has the fastest possible decreasing. In the subexponential case, 
$\hbox{\rm var}_n\psi \leq C \exp(-cn^{\gamma})$,
the optimal choice turns to be $n(r)=r^{1+\gamma}$, while in the polynomially decreasing 
case, $\hbox{\rm var}_n\psi\leq C n^{-q}$, the optimal choice is 
$n(r)=r^{(q-1)/(q-1-\epsilon)}$. This gives a bound in $n^{-q+2+\epsilon}$.
\end{remark}

\section{Comments and Open Questions}~\label{section-concluding-remarks}\

\medskip\noindent In our previous work~\cite{2003ChU} we made two restrictive assumptions, 
namely that $\psi$ is a locally constant potential and the image of the starting SFT under 
the amalgamation map $\pi$ is still a SFT (in general it is a sofic subshift).
In that setting, we could prove, under sufficient conditions, that 
$\mu_{\psi}\circ\pi^{-1}$ is a Gibbs measure for a H\"older continuous potential $\phi$. 
We also exhibited an example showing that one of our sufficient conditions turns out 
to be necessary in that otherwise the induced potential $\phi$ is not defined everywhere.

\noindent 
We conjecture the following: Let $\pi:A\to B$ be an amalgamation map as above, 
$X\subset A^{\mathbb N}$ a SFT and $Y\subset B^{\mathbb N}$ the resulting sofic subshift. 
Then the pushforward measure of a Gibbs measure for a H\"older continuous potential is 
a ``weak" Gibbs measure $\mu_\phi$ in that~\eqref{gibbs-inequality}
does not hold for every $\boldsymbol{a}$ but for almost all $\boldsymbol{a}$ (w.r.t. $\mu_\phi$).


\medskip\noindent\section{Proofs}

\medskip\noindent \subsection{Preliminary result: Birkhoff's refinement of 
Perron--Frobenius Theorem}~\label{refined-perron-frobenius}\

\medskip\noindent 
Let $E,E'$ be finite sets and  $M:E\times E'\to {\mathbb R}^+$ be a row allowable 
non--negative matrix, {\em i.e.}, a matrix such that $M\boldsymbol{x}>0$ whenever 
$\boldsymbol{x}>0$.  
Let us define the set
\[ \Delta_E:=\left\{\boldsymbol{x}\in \ ]0,1[^E:\ 
    |\boldsymbol{x}|_1:=\sum_{e\in E}\boldsymbol{x}(e)=1\right\}, \] 
and similarly $\Delta_{E'}$. We supply $\Delta_E$ with the distance
\[ \delta_E(x,y):=\max_{e,f\in E}
\log\frac{\boldsymbol{x}(e) \boldsymbol{y}(f)}{\boldsymbol{x}(f) \boldsymbol{y}(e)}. \] 
On $\Delta_{E'}$ we define $\delta_{E'}$ accordingly. Let us now define
\[ \tau(M):=\frac{1-\sqrt{\Phi(M)}}{1+\sqrt{\Phi(M)}} \]
where
\[ \Phi(M):=\left\{\begin{array}{ll} \displaystyle \min_{e,f\in E, \,e',f'\in E'}
              \frac{M(e,e')M(f,f')}{M(e,f')M(f,e')} 
                               & \text{ if } M > 0, \\ \\
                      0        & \text{otherwise.}
                \end{array}\right.\]
Here $M>0$ means that all entries of $M$ are strictly positive.

\medskip\noindent
\begin{theorem}[After Garrett Birkhoff]\label{teorema-contraccion}
Let $M:E\times E'\to {\mathbb R}^+$ be row allowable, and $F_M:\Delta_{E'}\to \Delta_E$ be 
such that
\[ F_M \boldsymbol{x}:=\frac{M\boldsymbol{x}}{|M\boldsymbol{x}|_1} 
          \hskip 10pt \text{ for each } \boldsymbol{x}\in \Delta_{E'}. \]
Then, for all $\boldsymbol{x},\boldsymbol{y}\in \Delta_{E'}$, we have
\[ \delta_{E}(F_M\boldsymbol{x},F_M\boldsymbol{y})\leq 
                    \tau(M) \delta_{E'}(\boldsymbol{x},\boldsymbol{y}). \]
We have $\tau(M)<1$ if and only if $M>0$. 
\end{theorem}

\medskip\noindent For a proof of this important result, see~\cite{caroll} for instance. 
It can also be deduced from the proof of a similar theorem concerning square matrices which 
can be found in~\cite{senetabook}. As a corollary of this result we obtain the following 
form of the Perron-Frobenius Theorem.

\begin{corollary}[Enhanced Perron--Frobenius Theorem]\label{teorema-perron-frobenius}
Suppose that $M:E\times E\to {\mathbb R}^+$ is primitive i.e., there exists 
$\ell \in {\mathbb N}$ such that $M^\ell>0$. Then its maximal eigenvalue $\rho_M$ is 
simple and it has a unique right eigenvector $\bar{R}_M\in \Delta_E$, and 
a unique left eigenvector $\bar{L}_M$ satisfying $\bar{L}_M^\dag \bar{R}_M=1$. 
Furthermore, for each $\boldsymbol{x}\in\Delta_E$ and each $n\in{\mathbb N}$ we have
\[M^n\boldsymbol{x} \lessgtr (\bar{L}_M^\dag \boldsymbol{x}) \rho_M^n  \bar{R}_M \, 
  \exp\left(
    \pm \frac{\ell\delta_{E}(\boldsymbol{x},F_M\boldsymbol{x})}{1-\tau}\ \tau^{\lfloor n/\ell\rfloor}
         \right),\]
with $\tau:=\tau(M^\ell)<1$.
\end{corollary}

\medskip\noindent \begin{proof} Let us first remark that $F_{M^\ell}=F^{\ell}_M$. Since 
$M^\ell>0$, then Theorem~\ref{teorema-contraccion} and the Contraction Mapping Theorem 
imply the existence of a unique fixed point $\boldsymbol{x}_M=F_M\boldsymbol{x}_M\in\Delta_E$ 
such that
\begin{eqnarray*}
\delta_E(F_{M}^n\boldsymbol{x},\boldsymbol{x}_M) & \leq & 
  \sum_{k=0}^\infty \delta_E
             \left(F_M^{n+k\ell}\boldsymbol{x},F_M^{n+(k+1)\ell}\boldsymbol{x}\right) \\ 
                                                 & \leq & 
  \frac{\delta_E(\boldsymbol{x},F_M^\ell\boldsymbol{x})\tau^{\lfloor n/\ell\rfloor}
       }{1-\tau} \leq
  \frac{\ell\delta_{E}(\boldsymbol{x},F_M\boldsymbol{x})\tau^{\lfloor n/\ell\rfloor}}{1-\tau},
\end{eqnarray*}
for each $n\in{\mathbb N}$ and $\boldsymbol{x}\in\Delta_E$. From the definition of 
projective distance it follows that, for each $\boldsymbol{x}\in\Delta_E$ and $n\in{\mathbb N}$ 
there is a constant $C=C(\boldsymbol{x},n)$ such that 
\begin{equation}\label{powers-convergence}
M^n\boldsymbol{x} \lessgtr C(\boldsymbol{x},n)\, \boldsymbol{x}_M \, 
\exp\left(\pm \frac{\ell\delta_{E}(\boldsymbol{x},F_M\boldsymbol{x})
                   }{1-\tau}\tau^{\lfloor n/\ell\rfloor} \right).
\end{equation}

\medskip\noindent Let us now prove that $\boldsymbol{x}_M\equiv\bar{R}_M\in\Delta_E$ is the 
unique positive right eigenvector associated to the maximum eigenvalue 
$\rho_M:=\max{\rm spec}(M)$. Indeed, since $F_M\boldsymbol{x}_M=\boldsymbol{x}_M$, then 
$M\boldsymbol{x}_M=\lambda\boldsymbol{x}_M$ for some $\lambda>0$. Now, if 
$M\boldsymbol{y}=\lambda\boldsymbol{y}$ for some $\boldsymbol{y}\in{\mathbb C}^E$, and taking 
into account that $M$ is a real matrix, then $\boldsymbol{y}=a\,\boldsymbol{z}$ for some 
$a\,\in {\mathbb C}$ and $\boldsymbol{z}\in\Delta_E$. Therefore $\lambda$ is a simple eigenvalue. 
It follows from Theorem~\ref{teorema-contraccion} and the contraction mapping theorem that 
$\boldsymbol{z}=\bar{R}_M$ is the associated eigenvector.

\noindent Consider the map 
$\boldsymbol{x}\mapsto \min_{e\in E}(M\boldsymbol{x})(e)/\boldsymbol{x}(e)$ on $\Delta_E$, 
and extend it to ${\rm clos}(\Delta_E)$ (the closure is taken with respect to the euclidean 
distance), by allowing values in the extended reals 
$\bar{{\mathbb R}}:={\mathbb R}\cup\{\infty\}$.~\footnote{Here we are following a standard 
argument which can be found in~\cite{senetabook} for instance.} 
The resulting transformation is upper semicontinuous, and therefore there exists 
$\boldsymbol{x}_0\in {\rm clos}(\Delta_E)$ attaining the supremum, {\it i.e.}, such that
\[
\rho:=\sup_{\boldsymbol{x}\in \Delta_E}\min_{e\in E}
       \frac{(M\boldsymbol{x})(e)}{\boldsymbol{x}(e)}=\min_{e\in E}\frac{(M\boldsymbol{x}_0)(e)}{\boldsymbol{x}_0(e)}.
\]
This supremum is an eigenvalue, and the point where it is attained is its corresponding 
positive eigenvector. Indeed, if $M\boldsymbol{x}_0\neq \rho\boldsymbol{x}_0$, {\em i.e.} 
if $(M\boldsymbol{x}_0)(e)>\rho\boldsymbol{x}_0(e)$ for some $e\in E$, then 
$M^{\ell+1}\boldsymbol{x}_0 > \rho\,M^\ell\boldsymbol{x}_0$ which implies that 
$\rho < \sup_{\boldsymbol{x}\in \Delta_E}\min_{e\in E}(M\boldsymbol{x})(e)/\boldsymbol{x}(e)$. 
Therefore $\boldsymbol{x}_0$ is a non--negative eigenvector for $M$, but since 
$M^\ell\boldsymbol{x}_0=\rho^\ell\boldsymbol{x}_0>0$, then necessarily 
$\boldsymbol{x}_0=\bar{R}_M$ and $\lambda=\rho$.

\noindent Finally, if $0\neq\boldsymbol{y}\in{\mathbb C}^E$ is a right eigenvector of $M$, 
associated to another eigenvalue $\lambda'\in {\mathbb C}$, then
\[
|\lambda'| \, |\boldsymbol{y}|=|M\boldsymbol{y}|\leq M|\boldsymbol{y}|,
\]
where $|\boldsymbol{z}|$ denotes the coordinatewise absolute value of the vector 
$\boldsymbol{z}\in {\mathbb C}^E$, and the inequality holds at each coordinate. If 
$|\lambda'| < \min_{e\in E}(|(M\boldsymbol{y})(e)|)/(|\boldsymbol{y}(e)|)$, 
we can find a vector $\boldsymbol{y}^+\in \Delta_E$ by slightly changing $|\boldsymbol{y}|$ 
at coordinates $e\in E$ where $\boldsymbol{y}(e)=0$ and then normalizing, so that 
$|\lambda'| \leq \min_{e\in E}(M\boldsymbol{y}^+)(e)/\boldsymbol{y}^+(e)$. 
If on the contrary $|\lambda'| = \min_{e\in E}|(M\boldsymbol{y})(e)|/|\boldsymbol{y}(e)|$, 
then $M^{\ell+1}|\boldsymbol{y}| \geq |\lambda'|M^\ell|\boldsymbol{y}|$, and normalizing 
$M^\ell|\boldsymbol{y}|$ we obtain $\boldsymbol{y}^+\in \Delta_E$ such a way that 
$|\lambda'| \leq \min_{e\in E}(M\boldsymbol{y}^+)(e)/\boldsymbol{y}^+(e)$. 
We conclude that,
\[|\lambda'| \leq \sup_{\boldsymbol{x}\in \Delta_E}
                \min_{e\in E}\frac{(M\boldsymbol{x})(e)}{\boldsymbol{x}(e)}:=\rho \]
for each $\lambda'\in{\rm spec}(M)$, therefore $\rho\equiv\rho_M:=\max{\rm spec}(M)$.

\bigskip\noindent It remains to prove that in Ineq.~\eqref{powers-convergence}, we have
$C(\boldsymbol{x},n)=(\bar{L}_M^\dag\boldsymbol{x})\,\rho_M^n$, where $\bar{L}_M>0$ 
is the left eigenvector associated to $\rho_M$, normalized so that $\bar{L}_M^\dag\bar{R}_M=1$. 
For this note that, by multiplying Ineq.~\eqref{powers-convergence} at left by $\bar{L}_M$, 
we obtain
\[ \rho_M^n(\bar{L}_M^\dag\boldsymbol{x})\lessgtr  (\bar{L}_M^\dag\bar{R}_M)\, 
       C(\boldsymbol{x},n)\, 
       \exp\left(\pm \frac{\ell\delta_{E}(\boldsymbol{x},F_M\boldsymbol{x})
                         }{1-\tau}\tau^{\lfloor n/\ell\rfloor} \right), \] 
hence $C(\boldsymbol{x},n)=\rho_M^n(\bar{L}_M^\dag\boldsymbol{x})/(\bar{L}_M^\dag\bar{R}_M)$, 
and the proof is finished.

\end{proof}

\bigskip\noindent \subsection{Proof of Lemma~\ref{lemma-projections}}~\label{proof-of-projections}\

\medskip\noindent

\medskip\noindent \subsubsection{The right eigenvector}\

\medskip\noindent Notice that the transition matrix ${\mathcal M}_r:={\mathcal M}_{\psi_r}$ 
is primitive with primitivity index $r$, hence, according to 
Corollary~\ref{teorema-perron-frobenius}
\[ {\mathcal M}_r^n\boldsymbol{x} \lessgtr (\bar{L}_r^\dag \boldsymbol{x})\rho_r^n \bar{R}_r \, 
\exp\left( \pm \frac{r\,\delta_{A^r}(\boldsymbol{x},F_r\boldsymbol{x})
                   }{1-\tau}\tau^{\lfloor\frac{n}{r}\rfloor}\right), \]
for each $\boldsymbol{x}\in\Delta_E$ and $n\in{\mathbb N}$. Here $\rho_r$ denotes the maximal
eigenvalue of ${\mathcal M}_r$, $\bar{R}_r\in \Delta_{A^r}$ its unique right eigenvector in 
the simplex, $\bar{L}_r$ its unique associated left eigenvector satisfying 
$\bar{L}_r^\dag \bar{R}_r=1$, and $\tau:=\tau({\mathcal M}_r^r)$ denotes the contraction 
coefficient associated to the positive matrix ${\mathcal M}_r^r$. 

\noindent Let us now obtain explicit an upper bound for $\tau$ and for the distance 
$\delta_{A^r}(\boldsymbol{x},F_r\boldsymbol{x})$ for particular values of 
$\boldsymbol{x}\in \Delta_E$. First,
\[ \Phi({\mathcal M}_r^r)\geq
  \min_{ \boldsymbol{u},\boldsymbol{v},\boldsymbol{u}',
         \boldsymbol{v}',\boldsymbol{u}'',\boldsymbol{v}''\in A^r} 
   \frac{ {\mathcal M}_{r}(\boldsymbol{u},\boldsymbol{u}')
          {\mathcal M}_{r}(\boldsymbol{v},\boldsymbol{v}'')  
       }{ {\mathcal M}_{r}(\boldsymbol{u},\boldsymbol{v}')
          {\mathcal M}_{r}(\boldsymbol{v},\boldsymbol{u}'')}
                          \geq \exp\left(-2\sum_{k=0}^r\hbox{\rm var}_k\psi\right) > 0. \]
Therefore $\tau \leq 1-e^{-\sum_{k=0}^r\hbox{\rm var}_k\psi}$ and 
$(1-\tau)^{-1}\leq e^{\sum_{k=0}^r\hbox{\rm var}_k\psi}$. 

\noindent 
Let $s_\psi:=\sum_{k=0}^\infty\hbox{\rm var}_k\psi$, and $\theta:=1-e^{-s_\psi}$. With this, 
and taking into account the upper bound for $\tau$ and $(1-\tau)^{-1}$, we obtain
\begin{equation}\label{rv-estimate}
{\mathcal M}_r^n\boldsymbol{x} \lessgtr (\bar{L}_r^\dag \boldsymbol{x}) \rho_r^n  \bar{R}_r \, 
       \exp\left( \pm r\,\delta_{A^r}(\boldsymbol{x},F_r\boldsymbol{x})\,
                 e^{s_\psi}\theta^{\left\lfloor \frac{n}{r}\right\rfloor} \right).
\end{equation}

\medskip\noindent 
On the other hand, for 
$\bar{u}:=(1/(\textup{Card}(A^r)),\ldots,1/(\textup{Card}(A^r)))^{\dag}\in \Delta_{A^r}$, 
we have
\begin{eqnarray*}
\delta_{A^r}(\bar{u},F_r\bar{u}) & := & 
     \max_{\boldsymbol{w},\boldsymbol{w}'\in A^r}
   \log\left( \frac{\bar{u}(\boldsymbol{w}')({\mathcal M}_r \bar{u})(\boldsymbol{w})
                  }{\bar{u}(\boldsymbol{w})({\mathcal M}_r\bar{u})(\boldsymbol{w}')}\right) \\
                                & \leq & r\,\log(\textup{Card}(A))+2\|\psi\| < r\,  C_0,
\end{eqnarray*}
where $C_0:=2\left(\log(\textup{Card}(A))+\|\psi\|\right)$. Therefore, by taking 
$\boldsymbol{x}=\bar{u}$ and $n=r^2$ in~\eqref{rv-estimate}, we finally obtain 
\begin{equation}\label{estimate-right}
\bar{R}_r(\boldsymbol{u})\lessgtr 
  \frac{\sum_{\boldsymbol{a}\in\hbox{\rm \tiny Per}_{r^2}(A^{\mathbb N})\cap [\boldsymbol{u}]}
            e^{S_{r^2-r-1}\psi_r(\boldsymbol{a})}
       }{\rho_r^{r^2}|\bar{L}_r|_1}\, e^{\pm C_0\,r^2\exp(s_\psi)\theta^r}.
\end{equation}

\medskip\noindent
\subsubsection{Ansatz for the induced potential}\

\medskip\noindent To each word $\boldsymbol{w}\in B^{r}$ we associate the simplex
\[\Delta_{\boldsymbol{w}}:=\left\{\boldsymbol{x}\in (0,1)^{E_{\boldsymbol{w}}}:\ 
       |\boldsymbol{x}|_1:=\sum_{\boldsymbol{v}\in E_{\boldsymbol{w}}}
         \boldsymbol{x}_{\boldsymbol{v}}=1\right\}, \]
where $E_{\boldsymbol{w}}:=\{\boldsymbol{v}\in A^r:\ \pi\boldsymbol{v}=\boldsymbol{w}\}$.

\noindent  
Let ${\mathcal M}_r,\rho_r,\bar{L}_r:=\bar{L}_{\psi_r}$ and $\bar{R}_r:=\bar{R}_{\psi_r}$ 
be as above, and define, for each $\boldsymbol{w}\in B^r$, the restriction 
$\bar{L}_{r,\boldsymbol{w}}:=\bar{L}_r|_{E_{\boldsymbol{w}}}\in (0,\infty)^{E_{\boldsymbol{w}}}$. 
Define $\bar{R}_{r,\boldsymbol{w}}$ in the analogous way, and for each $\boldsymbol{w}\in B^{r+1}$
let ${\mathcal M}_{r,\boldsymbol{w}}$ be the restriction of ${\mathcal M}_r$ to the coordinates 
in $E_{\boldsymbol{w}_0^{r-1}}\times E_{\boldsymbol{w}_1^r}$. Using this, and taking into account 
Eq.~\eqref{parry-measure}, which applies to our $(r+1)$--symbol potential $\psi_r$, 
we derive the matrix expression
\[ \nu_r[\boldsymbol{b}_0^n]\equiv 
 \sum_{\pi\boldsymbol{a}_0^n=\boldsymbol{b}_0^n} \mu_{\psi_r}[\boldsymbol{a}_0^n]
               =\bar{L}_{r,\boldsymbol{b}_0^{r-1}}^{\dag}
               \left(\frac{\prod_{j=0}^{n-r}{\mathcal M}_{r,\boldsymbol{b}_{j}^{j+r}}
                          }{\rho_r^{n-r+1}} \right)\bar{R}_{r,\boldsymbol{b}_{n-r+1}^n},
\]
for the induced measure $\nu_r:=\mu_{\psi_r}\circ\pi^{-1}$. It follows from this that
\begin{equation}\label{ansatz}
\log\left(\frac{\nu_r[\boldsymbol{b}_0^n]}{\nu_r[\boldsymbol{b}_1^n]}\right)=
\log\left(\frac{\left(\bar{L}_{r,\boldsymbol{b}_0^{r-1}}\right)^{\dag}
     \prod_{j=0}^{n-r}{\mathcal M}_{r,\boldsymbol{b}_{j}^{j+r}}\bar{R}_{r,\boldsymbol{b}_{n-r+1}^n}
              }{\left(\bar{L}_{r,\boldsymbol{b}_1^{r}}\right)^{\dag}
     \prod_{j=1}^{n-r}{\mathcal M}_{r,\boldsymbol{b}_{j}^{j+r}}\bar{R}_{r,\boldsymbol{b}_{n-r+1}^n}}
     \right)-\log(\rho_r).
\end{equation}

\medskip\noindent  For each $\boldsymbol{w}\in A^{r+s}$, with $s\geq 1$, let
${\mathcal M}_{r,\boldsymbol{w}}:=\prod_{j=0}^{s-1}{\mathcal M}_{r,\boldsymbol{w}_j^{j+r}}$, 
and define the transformation
$F_{r,\boldsymbol{w}}:\Delta_{\boldsymbol{w}_{s}^{r+s-1}}\to\Delta_{\boldsymbol{w}_0^{r-1}}$
such that
\[F_{r,\boldsymbol{w}}\boldsymbol{x}=
      \frac{{\mathcal M}_{r,\boldsymbol{w}}\boldsymbol{x}
          }{|{\mathcal M}_{r,\boldsymbol{w}}\boldsymbol{x}|_1}.\]
For each $\boldsymbol{b}\in B^{\mathbb N}$ and $s,t\in{\mathbb N}$, let
\[\boldsymbol{x}_{r,\boldsymbol{b}_{s+1}^{s+t+r}}:=
  F_{r,\boldsymbol{b}_{s+1}^{s+r+1}}\circ\cdots\circ
  F_{r,\boldsymbol{b}_{s+t}^{s+t+r}}
  \left(\bar{R}_{r,\boldsymbol{b}_{s+t+1}^{s+t+r}}/
          |\bar{R}_{r,\boldsymbol{b}_{s+t+1}^{s+t+r}}|_1\right)\in
                                            \Delta_{\boldsymbol{b}_{s+1}^{s+r}}. \]
By convention, 
$\boldsymbol{x}_{r,\boldsymbol{b}_{s+1}^{s+r}}\equiv
\bar{R}_{r,\boldsymbol{b}_{s+1}^{s+r}}/|\bar{R}_{r,\boldsymbol{b}_{s+1}^{s+r}}|_1\in 
\Delta_{\boldsymbol{b}_{s+1}^{s+r}}$. 
Using this notation, and after the adequate renormalization, Eq.~\eqref{ansatz} becomes
\begin{equation}\label{ansatz-2}
\log\left(\frac{\nu_r[\boldsymbol{b}_0^n]}{\nu_r[\boldsymbol{b}_1^n]}\right)= 
\log\left(
        \frac{\left(\bar{L}_{r,\boldsymbol{b}_0^{r-1}}\right)^{\dag}
               {\mathcal M}_{r,\boldsymbol{b}_0^r}\boldsymbol{x}_{r,\boldsymbol{b}_{1}^n} 
            }{\left(\bar{L}_{r,\boldsymbol{b}_1^{r}}\right)^{\dag}
                  \boldsymbol{x}_{r,\boldsymbol{b}_{1}^n}} \right)-\log(\rho_r).
\end{equation}

\medskip\noindent
\subsubsection{Convergence of the inhomogeneous product}\

\medskip\noindent Let us now prove the convergence of the sequence
$\left(\boldsymbol{x}_{\boldsymbol{b}_1^n}\right)_{n\geq r}$. For this notice that
\begin{eqnarray*}
\boldsymbol{x}_{r,\boldsymbol{b}_1^n}&:=&
F_{r,\boldsymbol{b}_{1}^{r+1}}\circ\cdots\circ F_{r,\boldsymbol{b}_{n-r}^{n}}
                                 \boldsymbol{x}_{r,\boldsymbol{b}_{n-r+1}^{n}}\\
                                      &=&
F_{r,\boldsymbol{b}_1^n}\boldsymbol{x}_{r,\boldsymbol{b}_{n-r+1}^n}
                                       =
F_{r,\boldsymbol{b}_{1}^{2r}}\circ F_{r,\boldsymbol{b}_{r+1}^{3r}}\circ\cdots\circ 
F_{r,\boldsymbol{b}_{(k-1)r+1}^{(k+1)r}}\boldsymbol{x}_{r,\boldsymbol{b}_{kr+1}^n},
\end{eqnarray*}
where $k:=\left\lfloor \frac{n}{r}\right\rfloor-1$. Now, since 
${\mathcal M}_{r,\boldsymbol{w}}>0$ for each $\boldsymbol{w}\in B^{2r}$, then 
Theorem~\ref{teorema-contraccion} ensures that the associated transformation
$F_{r,\boldsymbol{w}}:\Delta_{\boldsymbol{w}_r^{2r-1}}\to \Delta_{\boldsymbol{w}_0^{r-1}}$,
is a contraction with contraction coefficient
$\tau_{\boldsymbol{w}}=(1-\sqrt{\Phi_{\boldsymbol{w}}})/(1+\sqrt{\Phi_{\boldsymbol{w}}})$,
where
\begin{eqnarray}\label{Phi-omega}
\Phi_{\boldsymbol{w}}&:=  & 
\min_{\boldsymbol{v},\boldsymbol{u} \in E_{\boldsymbol{w}_0^{r-1}},\,
      \boldsymbol{v}',\boldsymbol{u}'\in E_{\boldsymbol{w}_r^{2r-1}}}
      \frac{ {\mathcal M}_{r,\boldsymbol{w}}(\boldsymbol{v},\boldsymbol{v}')
             {\mathcal M}_{r,\boldsymbol{w}}(\boldsymbol{u},\boldsymbol{u}')
          }{ {\mathcal M}_{r,\boldsymbol{w}}(\boldsymbol{v},\boldsymbol{u}')
             {\mathcal M}_{r,\boldsymbol{w}}(\boldsymbol{u},\boldsymbol{v}')} \nonumber\\
                    &\geq & 
\min_{\boldsymbol{v},\boldsymbol{u}, \boldsymbol{v}',
      \boldsymbol{u}', \boldsymbol{v}'',\boldsymbol{u}''\in A^r}
     \frac{{\mathcal M}_{r,\boldsymbol{w}}(\boldsymbol{v},\boldsymbol{v}')
           {\mathcal M}_{r,\boldsymbol{w}}(\boldsymbol{u},\boldsymbol{u}'') 
          }{{\mathcal M}_{r,\boldsymbol{w}}(\boldsymbol{v},\boldsymbol{u}')
            {\mathcal M}_{r,\boldsymbol{w}}(\boldsymbol{u},\boldsymbol{v}'')} \nonumber\\
         &\geq &\exp\left(-2\sum_{k=0}^r\hbox{\rm var}_k\psi\right) \geq e^{-2\,s_\psi}> 0.
\end{eqnarray}
Recall that $s_\psi=\sum_{k=0}^\infty \hbox{\rm var}_k\psi$. From Ineq.~\eqref{Phi-omega} we 
obtain a uniform upper bound for the contraction coefficients, 
$\tau_{\boldsymbol{w}}\leq \theta:=1-\exp(-s_\psi) < 1$, which allows us to
establish the uniform convergence of the sequence 
$\left(\boldsymbol{x}_{\boldsymbol{b}_1^n}\right)_{n\geq r}$ with respect to 
$\boldsymbol{b}\in B^{\mathbb N}$. Indeed, for $\boldsymbol{b}\in B^{\mathbb N}$ fixed and 
$m > n$, we have
\begin{equation}\label{convergence-1}
\delta_{E_{\boldsymbol{b}_1^r}}\left(\boldsymbol{x}_{\boldsymbol{b}_1^n},
                                     \boldsymbol{x}_{\boldsymbol{b}_1^m} \right)
   \leq \theta^k\, \delta_{E_{\boldsymbol{b}_{kr+1}^{(k+1)r}}} 
                      \left(\boldsymbol{x}_{r,\boldsymbol{b}_{kr+1}^n},
                      \boldsymbol{x}_{r,\boldsymbol{b}_{kr+1}^{m}}\right)
\end{equation}
where $k:=\lfloor \frac{n}{r}\rfloor-1$. On the other hand,
\begin{eqnarray*}
\delta_{E_{\boldsymbol{b}_{kr+1}^{(k+1)r}}}
       \left(\boldsymbol{x}_{r,\boldsymbol{b}_{kr+1}^n},
             \boldsymbol{x}_{r,\boldsymbol{b}_{kr+1}^{m}}\right) &\leq & \sum_{j=0}^{k'}
\delta_{E_{\boldsymbol{b}_{kr+1}^{(k+1)r}}}
      \left(\boldsymbol{x}_{r,\boldsymbol{b}_{kr+1}^{n+jr}},
            \boldsymbol{x}_{r,\boldsymbol{b}_{kr+1}^{n+(j+1)r}}\right)\\
                                                                  &     & +
\delta_{E_{\boldsymbol{b}_{kr+1}^{(k+1)r}}}
      \left(\boldsymbol{x}_{r,\boldsymbol{b}_{kr+1}^{n+(k'+1)r}},
            \boldsymbol{x}_{r,\boldsymbol{b}_{kr+1}^{m}}\right),
\end{eqnarray*}
where $k':=\lfloor (m-n)/r\rfloor-1$. By convention, when $k'=-1$, the summation 
in the right--hand side is zero. Then, since all the matrices
${\mathcal M}_{r,\boldsymbol{w}}$ are row allowable and positive for $\boldsymbol{w}\in B^{2r}$, 
then we have
\begin{equation}\label{convergence-2}
\delta_{E_{\boldsymbol{b}_{kr+1}^{(k+1)r}}}
            \left(\boldsymbol{x}_{r,\boldsymbol{b}_{kr+1}^n},
                  \boldsymbol{x}_{r,\boldsymbol{b}_{kr+1}^{m}}\right) \leq T_1+T_2+T_3
\end{equation}
where
\begin{equation}\label{T1}
T_1:=
\delta_{E_{\boldsymbol{b}_{n-r+1}^{n}}}
  \left(\boldsymbol{x}_{r,\boldsymbol{b}_{n-r+1}^{n}},
        F_{r,\boldsymbol{b}_{n-r+1}^{n+r}}\boldsymbol{x}_{r,\boldsymbol{b}_{n+1}^{n+r}}\right),
\end{equation}
\begin{equation}\label{T2}
T_2:=
\sum_{j=1}^{k'} \, \theta^j\,
\delta_{E_{\boldsymbol{b}_{n+(j-1)r+1}^{n+jr}}}
\left(\boldsymbol{x}_{r,\boldsymbol{b}_{n+(j-1)r+1}^{n+jr}},
 F_{r,\boldsymbol{b}_{n+(j-1)r+1}^{n+(j+1)r}}\boldsymbol{x}_{r,\boldsymbol{b}_{n+jr+1}^{n+(j+1)r}}
\right)
\end{equation}
and 
\begin{equation}\label{T3}
T_3:=
\theta^{k'}\,\delta_{E_{\boldsymbol{b}_{n+k'r+1}^{n+(k'+1)r}}}
       \left(\boldsymbol{x}_{r,\boldsymbol{b}_{n+k'r+1}^{n+(k'+1)r}},
       F_{r,\boldsymbol{b}_{n+k'r+1}^{m}}\boldsymbol{x}_{r,\boldsymbol{b}_{m-r+1}^{m}}\right).
\end{equation}
Once again, we convene that $T_2=0$ if $k'=-1$.

\medskip\noindent Now, for each $\boldsymbol{w},\boldsymbol{w}'\in B^r$, and 
$\boldsymbol{v}\in B^s$ with $r< s< 2r$, and such that 
$\boldsymbol{v}_0^{r-1}=\boldsymbol{w},\,\boldsymbol{v}_{s-r+1}^s=\boldsymbol{w}'$, we have
\begin{eqnarray*}
\delta_{\boldsymbol{w}}(\boldsymbol{x}_{r,\boldsymbol{w}},
                        F_{r,\boldsymbol{v}}\boldsymbol{x}_{r,\boldsymbol{w}'}) 
                                  & = & 
\max_{\boldsymbol{u},\boldsymbol{u}'\in E_{\boldsymbol{w}}}
\log\left(\frac{
    \boldsymbol{x}_{r,\boldsymbol{w}}(\boldsymbol{u})\ 
             (F_{r,\boldsymbol{v}}\boldsymbol{x}_{r,\boldsymbol{w}'})(\boldsymbol{u}')
              }{
   \boldsymbol{x}_{r,\boldsymbol{w}}(\boldsymbol{u}')\ 
             (F_{r,\boldsymbol{v}}\boldsymbol{x}_{r,\boldsymbol{w}'})(\boldsymbol{u})}\right)\\
                                 & \leq & 
\max_{\boldsymbol{u},\boldsymbol{u}'\in E_{\boldsymbol{w}}}
\log\left(\frac{
     \bar{R}_{r}(\boldsymbol{u})\ 
        ({\mathcal M}_{r,\boldsymbol{v}}\bar{R}_{r,\boldsymbol{w}'})(\boldsymbol{u}')
               }{
     \bar{R}_{r}(\boldsymbol{u}')\ 
         ({\mathcal M}_{r,\boldsymbol{v}}\bar{R}_{r,\boldsymbol{w}'})(\boldsymbol{u})}\right).
\end{eqnarray*}
Hence, using the estimate for the right eigenvectors given in
Eq.~\eqref{estimate-right}, it follows that
\begin{eqnarray*}
\delta_{\boldsymbol{w}}(
 \boldsymbol{x}_{r,\boldsymbol{w}},F_{r,\boldsymbol{v}}\boldsymbol{x}_{r,\boldsymbol{w}'})
                                   &\leq& 
\max_{\boldsymbol{u},\boldsymbol{u}'\in A^r} 
\log \left(\frac{\displaystyle \sum_{\boldsymbol{a}\in\hbox{\rm \tiny Per}_{r^2}(A^{\mathbb N})}
       e^{S_{r^2-r-1}\psi_r(\boldsymbol{a})}
       \hskip -15pt \sum_{\boldsymbol{a}\in\hbox{\rm \tiny Per}_{r^2+s-r}(A^{\mathbb N})}
       e^{S_{r^2+s-2r-1}\psi_r(\boldsymbol{a})}
               }{\displaystyle \min_{\boldsymbol{a}\in\hbox{\rm \tiny Per}_{r^2}(A^{\mathbb N})}
       e^{S_{r^2-r-1}\psi_r(\boldsymbol{a})}\hskip -15pt
             \min_{\boldsymbol{a}\in\hbox{\rm \tiny Per}_{r^2+s-r}(A^{\mathbb N})}
       e^{S_{r^2+s-2r-1}\psi_r(\boldsymbol{a})}}\right)\\
                                    & + &
     2\,r^2\,C_0\,e^{s_\psi}\,\theta^r \leq  2\,r(r+1)\,C_0(e^{s_\psi}\,\theta^r+1),
\end{eqnarray*}
with $C_0=2(\log(\textup{Card}(A))+\|\psi\|)$ and $\theta=1-\exp(-s_\psi)$ as in 
Eq.~\eqref{estimate-right}. Using this upper bound in~\eqref{T1},~\eqref{T2} and~\eqref{T3}, 
we obtain from~\eqref{convergence-2}
\[
\delta_{E_{\boldsymbol{b}_{kr+1}^{(k+1)r}}}
\left(\boldsymbol{x}_{r,\boldsymbol{b}_{kr+1}^n},
      \boldsymbol{x}_{r,\boldsymbol{b}_{kr+1}^{m}}\right)
    \leq 2\,r(r+1)\,C_0(e^{s_\psi}\,\theta^r+1)\left(\theta^k+\frac{1}{1-\theta}\right),
\]
and with this, Ineq.~\eqref{convergence-1} becomes 
\begin{equation}\label{convergence-3}
\delta_{E_{\boldsymbol{b}_1^r}}
\left(\boldsymbol{x}_{\boldsymbol{b}_1^n},
       \boldsymbol{x}_{\boldsymbol{b}_1^m} \right)\leq 
        2\,r(r+1)\,C_0 (e^{s_\psi}\,\theta^r+1)\theta^{\lfloor \frac{n}{r}\rfloor-1} 
         \left(\theta^{\lfloor \frac{n}{r}\rfloor-1}+\frac{1}{1-\theta}\right)
\end{equation}
which holds for all $\boldsymbol{b}\in B^{\mathbb N}$ and $r <n < m$. Hence, 
$(\boldsymbol{x}_{\boldsymbol{b}_1^n})_{n\geq r}$ is a Cauchy sequence in complete space 
$\Delta_{\boldsymbol{b}_{1}^r}$, and the existence of the limit 
$\boldsymbol{x}_{\boldsymbol{b}_1^{\infty}}:=\lim_{m\to\infty}\boldsymbol{x}_{\boldsymbol{b}_1^m}$ 
is ensured for each $\boldsymbol{b}\in B^{\mathbb N}$. Furthermore, from Eq.~\eqref{convergence-3} 
it follows that
\begin{equation}\label{convergence-final}
\delta_{E_{\boldsymbol{b}_1^r}}
     \left(\boldsymbol{x}_{\boldsymbol{b}_1^n},
           \boldsymbol{x}_{\boldsymbol{b}_1^{\infty}} \right)\leq 
       2\,r(r+1)\,C_0(e^{s_\psi}\,\theta^r+1)\theta^{\lfloor \frac{n}{r}\rfloor-1}
         \left(\theta^{\lfloor \frac{n}{r}\rfloor-1}+\frac{1}{1-\theta}\right)
      \leq C_1\,r^2\,\theta^{\frac{n}{r}},
\end{equation}
with $C_1:=4C_0(1+e^{s_\psi}\,\theta)/(\theta^2(1-\theta))$.

\medskip\noindent
\subsubsection{The induced potential and the Gibbs inequality}\

\medskip\noindent Taking the Eq.~\eqref{convergence-final}, it follows that the limit
\begin{equation}\label{expression-potential}
\phi_{r}(\boldsymbol{b})=\lim_{n\to\infty}
\log\left(\frac{\nu_r[\boldsymbol{b}_0^n]}{\nu_r[\boldsymbol{b}_1^n]}\right)
                        =
\log\left(\frac{ 
 \left(\bar{L}_{r,\boldsymbol{b}_0^{r-1}}\right)^{\dag}{\mathcal M}_{r,\boldsymbol{b}_0^r}
          \boldsymbol{x}_{r,\boldsymbol{b}_{1}^\infty} 
              }{\left(\bar{L}_{r,\boldsymbol{b}_1^{r}}\right)^{\dag}
          \boldsymbol{x}_{r,\boldsymbol{b}_{1}^\infty}}\right) -\log(\rho_r),
\end{equation}
exists for each $\boldsymbol{b}\in B^{\mathbb N}$, and defines a continuous function 
$\boldsymbol{b}\mapsto\phi_{r}(\boldsymbol{b})$. This proves that the limit \eqref{lephir} 
in the statement of the lemma does exist. It remains to find an upperbound to its modulus 
of continuity.

\noindent Inequality~\eqref{convergence-final}, and the fact that 
$|\boldsymbol{x}_{\boldsymbol{b}_1^n}|_1=|\boldsymbol{x}_{\boldsymbol{b}_1^{\infty}}|_1=1$, 
imply that 
\[
\boldsymbol{x}_{\boldsymbol{b}_1^n}\lessgtr \boldsymbol{x}_{\boldsymbol{b}_1^{\infty}}\,
\exp\left(\pm\,C_1 r^2\theta^{\frac{n}{r}}\right)
\] 
for all $\boldsymbol{b}\in B^{\mathbb N}$ and $n > r$. With this, and taking into account 
Eqs.~\eqref{ansatz-2}
and~\eqref{expression-potential}, it follows that
\begin{eqnarray*}
\left| \phi_{r}(\boldsymbol{b})-
  \log\left(\frac{\nu_r[\boldsymbol{b}_0^n]}{\nu_r[\boldsymbol{b}_1^n]}\right)\right| 
                                               &\leq & 
\left| \log\left( 
           \frac{\left(\bar{L}_{r,\boldsymbol{b}_0^{r-1}}\right)^{\dag} 
                {\mathcal M}_{r,\boldsymbol{b}_0^r}\boldsymbol{x}_{r,\boldsymbol{b}_{1}^\infty} 
               }{\left(\bar{L}_{r,\boldsymbol{b}_0^{r}}\right)^{\dag}
                 {\mathcal M}_{r,\boldsymbol{b}_0^r}\boldsymbol{x}_{r,\boldsymbol{b}_{1}^n} 
                 }\right)- 
       \log\left(
          \frac{\left(\bar{L}_{r,\boldsymbol{b}_1^{r-1}}\right)^{\dag} 
                  \boldsymbol{x}_{r,\boldsymbol{b}_{1}^n} 
              }{ \left(\bar{L}_{r,\boldsymbol{b}_1^{r}}\right)^{\dag}
                  \boldsymbol{x}_{r,\boldsymbol{b}_{1}^\infty}}\right) \right|\\ \\
                                              &\leq&  C\,r^2\theta^{\frac{n}{r}},
\end{eqnarray*}
for all $\boldsymbol{b}\in B^{\mathbb N}$, $n >r $ and
$C:=2C_1=8C_0(1+e^{s_\psi}\,\theta)/(\theta^2(1-\theta))$. 
This proves~\eqref{regularitephir} in the statement of the lemma.

\noindent 
From this it can be easily deduced that $\nu_r\equiv\mu_{\psi_r}\circ\pi^{-1}$ satisfies 
the Gibbs Inequality~\eqref{gibbs-inequality} with potential $\phi_{r}$ and constants 
$P(\phi_{r},B^{\mathbb N})=0$ and
\[
C(\phi_{r},B^{\mathbb N})=\max_{\boldsymbol{b}\in B^{\mathbb N}}
\left(
\frac{\exp(S_{r^2}\phi_{r}(\boldsymbol{b}))}{\nu_r\left[\boldsymbol{b}_0^{r^2}\right]},
\frac{\nu_r\left[\boldsymbol{b}_0^{r^2}\right]}{\exp(S_{r^2}\phi_{r}(\boldsymbol{b}))}
\right)\ \exp\left(\frac{C\,r^2\,\theta^r}{1-\theta^{\frac{1}{r}}}\right).
\]
This proves the first statement of the lemma the proof of which is now complete.
\hfill{$\Box$}

\medskip\noindent \begin{remark}
As mentioned above (see \eqref{topological-pressure}), the topological pressure of
$\psi$ is given by
\[ P(\psi)=P(\psi,A^{\mathbb N})=\lim_{n\to\infty}\frac{1}{n}\log 
                \left(\sum_{\boldsymbol{a}\in\hbox{\rm \tiny Per}_{n}(A^{\mathbb N})}
                e^{S_n\psi(\boldsymbol{a})} \right). \]
Since $\psi\lessgtr \psi_r\pm \hbox{\rm var}_r\psi$, we get
\begin{eqnarray*}
\log(\rho_r) &     =  & 
\lim_{n\to\infty}\frac{1}{n}\log\left({\rm Trace}\left({\mathcal M}_r^n\right)\right)\\
             &\lessgtr& 
\lim_{n\to\infty}\frac{1}{n}\log\left(
                     \sum_{\boldsymbol{a}\in\hbox{\rm \tiny Per}_{n}(A^{\mathbb N})}
                        e^{S_n\psi(\boldsymbol{a})}\right)\pm \hbox{\rm var}_r\psi \\
              &\lessgtr& P(\psi)\pm \hbox{\rm var}_r\psi,
\end{eqnarray*}
for each $r\in{\mathbb N}$. 
\end{remark}

\bigskip\noindent 
\subsection{Proof of Lemma~\ref{lemma-projective}}~\label{proof-of-projective}\

\medskip\noindent
\subsubsection{Periodic approximations}\

\medskip\noindent Each Markov approximant $\mu_{\psi_r}$ can be seen as the limit of 
measures supported on periodic points as follows. Fix $n,r\in {\mathbb N}$ with 
$n\geq r$, and $\boldsymbol{w}\in A^n$. Then, for each $p > r+n$ we have,
\[{\mathcal P}_r^{(p)}[\boldsymbol{w}]:=\frac{
\sum_{\boldsymbol{a}\in\hbox{\rm \tiny Per}_p(A^{\mathbb N})\cap\,[\boldsymbol{w}]}
                                              e^{S_p\psi_r(\boldsymbol{a})}
                                            }{ 
\sum_{\boldsymbol{a}\in\hbox{\rm \tiny Per}_p(A^{\mathbb N})}e^{S_p\psi_r(\boldsymbol{a})}}.
\]
We can rewrite the above equation as
\begin{eqnarray*}
{\mathcal P}_r^{(p)}[\boldsymbol{w}]&=& 
\frac{\left(\prod_{s=0}^{n-r-1}
    {\mathcal M}_r(\boldsymbol{w}_{s}^{s+r-1},\boldsymbol{w}_{s+1}^{s+r})\right)
    \bar{e}_{\boldsymbol{w}_{n-r}^{n-1}}^{\dag}{\mathcal M}_r^{p-n+r}
    \bar{e}_{\boldsymbol{w}_0^{r-1}} 
    }{ \sum_{\zeta\in A^r} \bar{e}_{\zeta}^{\dag}{\mathcal M}_r^{p}\bar{e}_{\zeta}} \\
                                    &=&
\frac{\left(
  \prod_{s=0}^{n-r-1}{\mathcal M}_r(\boldsymbol{w}_{s}^{s+r-1},\boldsymbol{w}_{s+1}^{s+r})
      \right)\bar{e}_{\boldsymbol{w}_{n-r}^{n-1}}^{\dag}{\mathcal M}_r^{p-n}
  \left({\mathcal M}_r^{r}\bar{e}_{\boldsymbol{w}_0^{r-1}}\right) 
    }{ \sum_{\zeta\in A^r} \bar{e}_{\zeta}^{\dag} {\mathcal M}_r^{p-n} 
                   \left({\mathcal M}_r^{n}\bar{e}_{\zeta}\right)},
\end{eqnarray*}
with ${\mathcal M}_r$ as above, and $\bar{e}_{\zeta}\in \{0,1\}^{A^r}$ the vector with $1$ 
at coordinate $\zeta$ and zeros everywhere else. Now, since 
${\mathcal M}_r^{k}\bar{e}_{\zeta} > 0$ for each $k\geq r$ 
and $\zeta\in A^r$, then Corollary~\ref{teorema-perron-frobenius} applies, and 
using~\eqref{parry-measure} we obtain
\begin{eqnarray*}
{\mathcal P}_r^{(p)}[\boldsymbol{w}] &\lessgtr& 
\frac{ \bar{L}_r^{\dag}\left({\mathcal M}_r^{r}\bar{e}_{\boldsymbol{w}_0^{r-1}}\right)\, 
   \left(
   \prod_{s=0}^{n-r-1}{\mathcal M}_r(\boldsymbol{w}_{s}^{s+r-1},\boldsymbol{w}_{s+1}^{s+r})
   \right)\bar{e}_{\boldsymbol{w}_{n-r}^{n-1}}^{\dag}\bar{R}_r 
   }{\sum_{\zeta\in A^r} \bar{L}_r^{\dag} 
   \left({\mathcal M}_r^{n} \bar{e}_{\boldsymbol{w}'}\right)\,\bar{e}_{\zeta}^{\dag} \bar{R}_r
    }\ e^{\pm \frac{r\,D_0}{1-\tau}\,\tau^{\frac{p-n}{r}-2}}\\
                                     &\lessgtr&
\frac{\bar{L}_r^{\dag}(\boldsymbol{w}_0^{r-1})\, 
  \left(
  \prod_{s=0}^{n-r-1}{\mathcal M}_r(\boldsymbol{w}_{s}^{s+r-1},\boldsymbol{w}_{s+1}^{s+r})
  \right)\bar{R}_r(\boldsymbol{w}_{n-r}^{n-1}) 
    }{ \rho_r^{n-r}}\ e^{\pm \frac{r\,D_0}{1-\tau}\, \tau^{\frac{p-n}{r}-2}} \\
                                     &\lessgtr&
\mu_{\psi_r}[\boldsymbol{w}]\,\exp\left(\pm\frac{r\,D_0}{1-\tau}\,\tau^{\frac{p-n}{r}-2}\right),                  
\end{eqnarray*}
with $\bar{R}_r, \bar{L}_r, \rho_r$ and $\tau:=\tau({\mathcal M}_r^r)$ as before, and 
\[D_0:=2\,\max_{\zeta\in A^r}\delta_{A^r}
\left({\mathcal M}_r^r\bar{e}_{\zeta},{\mathcal M}_r^{r+1}\bar{e}_{\zeta}\right).\]

\medskip\noindent Since $\tau=(1-\sqrt{\Phi})/(1+\sqrt{\Phi})$, with
\[\Phi:= \min_{\boldsymbol{v},\boldsymbol{u},\, \boldsymbol{v}',\boldsymbol{u}'\in A^r} 
\frac{{\mathcal M}_{r,\boldsymbol{w}}(\boldsymbol{v},\boldsymbol{v}')
      {\mathcal M}_{r,\boldsymbol{w}}(\boldsymbol{u},\boldsymbol{u}') 
     }{{\mathcal M}_{r,\boldsymbol{w}}(\boldsymbol{v},\boldsymbol{u}')
       {\mathcal M}_{r,\boldsymbol{w}}(\boldsymbol{u},\boldsymbol{v}')}  
       \geq \exp\left(-2\sum_{k=0}^{r}\hbox{\rm var}_k\psi\right) > 0,\]
then $\tau \leq \theta:=1-\exp\left(-s_\psi\right)$, and
$(1-\tau)^{-1} \leq \exp\left(s_\psi\right)$, with 
$s_\psi:=\sum_{k=0}^{\infty}\hbox{\rm var}_k\psi$ as in Lemma~\ref{lemma-projections}. 
On the other hand, 
\begin{eqnarray*}
D_0 &\leq& 
2 \max_{\boldsymbol{w}',\boldsymbol{u},\boldsymbol{u}'\in A^r} 
\log\left(
   \frac{{\mathcal M}_{r}^r(\boldsymbol{u},\boldsymbol{w}')
         {\mathcal M}_{r}^{r+1}(\boldsymbol{u}',\boldsymbol{w}')
        }{{\mathcal M}_{r}^{r+1}(\boldsymbol{u},\boldsymbol{w}')
          {\mathcal M}_{r}^r(\boldsymbol{u}',\boldsymbol{w}')}\right)\\
    &\leq& 4\left(\log(\textup{Card}(A))+\sum_{k=1}^{r}\hbox{\rm var}_k\psi+\|\psi\|\right)\\
    & \leq & 4 \left(\log(\textup{Card}(A))+s_\psi+\|\psi\|\right)=:D_1.
\end{eqnarray*}
Using this explicit bound just obtained, we deduce the inequalities
\[\mu_{\psi_r}[\boldsymbol{w}]\lessgtr 
\frac{\sum_{\boldsymbol{a}\in\hbox{\rm \tiny Per}_p(A^{\mathbb N})\cap\,[\boldsymbol{w}]} 
  e^{S_p\psi_r(\boldsymbol{a})}
    }{\sum_{\boldsymbol{a}\in\hbox{\rm \tiny Per}_p(A^{\mathbb N})}
  e^{S_p\psi_r(\boldsymbol{a})}}\ \exp\left(\pm D_1\,r\,e^{s_\psi}\theta^{\frac{p-n}{r}-2}\right)
\]
for each $\boldsymbol{w}\in A^n$, and all $p > n+r$. It is easy to check that these inequalities
extend to each $\boldsymbol{w}\in \cup_{k=1}^nA^k$, and we finally obtain
\begin{equation}\label{periodic-approximation2}
\mu_{\psi_r}[\boldsymbol{w}]\lessgtr 
\frac{\sum_{\boldsymbol{a}\in\hbox{\rm \tiny Per}_p(A^{\mathbb N})\cap\,[\boldsymbol{w}]} 
e^{S_p\psi_r(\boldsymbol{a})}
     }{\sum_{\boldsymbol{a}\in\hbox{\rm \tiny Per}_p(A^{\mathbb N})}e^{S_p\psi_r(\boldsymbol{a})}}
                \ \exp\left(\pm D_1\,r\,e^{s_\psi}\theta^{\frac{p-\max(n,r)}{r}-2}\right)
\end{equation}
for all $r,n\in{\mathbb N}$, and $\boldsymbol{w}\in A^n$.

\medskip\noindent
\subsubsection{Telescopic product}\

\medskip\noindent Let us now compare two consecutive Markov approximants. Fix $n,r > 0$, 
and $\boldsymbol{w}\in A^n$. Then, for each $p > n+r+1$, 
Inequalities~\eqref{periodic-approximation2} ensure that
\[\frac{\mu_{\psi_r}[\boldsymbol{w}]}{\mu_{\psi_{r+1}}[\boldsymbol{w}]}
\lessgtr 
\frac{\sum_{\boldsymbol{a}\in\hbox{\rm \tiny Per}_p(A^{\mathbb N})\cap\,[\boldsymbol{w}]} 
     e^{S_p\psi_r(\boldsymbol{a})} 
     }{\sum_{\boldsymbol{a}\in\hbox{\rm \tiny Per}_p(A^{\mathbb N})\cap\,[\boldsymbol{w}]}
     e^{S_p\psi_{r+1}(\boldsymbol{a})}}\,
\frac{\sum_{\boldsymbol{a}\in\hbox{\rm \tiny Per}_p(A^{\mathbb N})}
     e^{S_p\psi_{r+1}(\boldsymbol{a})}
    }{\sum_{\boldsymbol{a}\in\hbox{\rm \tiny Per}_p(A^{\mathbb N})}
     e^{S_p\psi_{r}(\boldsymbol{a}) }}\,\exp\left(\pm C\, r\,\theta^{\frac{p-q}{r}-2}\right),
\]
with $q=\max(r+1,n)$ and $C:=2\,e^{s_\psi}\, D_1$. 
Since $\psi_{r+1}\lessgtr \psi_{r}\pm \hbox{\rm var}_{r+1}\psi$, then we have
\begin{eqnarray*} 
\frac{\mu_{\psi_r}[\boldsymbol{w}]}{\mu_{\psi_{r+1}}[\boldsymbol{w}]}
&\lessgtr&
\exp\left(\pm\left(2\,p\,\hbox{\rm var}_{r+1}\psi+C\,r\,\theta^{\frac{p-q}{r}-2}\right)\right)\\
&\lessgtr&
\exp\left(\pm \left(2\,p\,\hbox{\rm var}_{r+1}\psi+C\,r\,\theta^{\frac{p-r-n-1}{r}-2}\right)\right)
\end{eqnarray*}
for all $r\in{\mathbb N}$, $\boldsymbol{w}\in \cup_{k=1}^nA^k$, and $p > n+r+1$. 
Let $p=(r+1)(r+2)+n-1$, then for each 
$r' > r\in {\mathbb N}$ and $\boldsymbol{w}\in \cup_{k=1}^nA^k$  we have
\[\frac{\mu_{\psi_r}[\boldsymbol{w}]}{\mu_{\psi_{r'}}[\boldsymbol{w}]}
\lessgtr 
\exp\left(
\pm\,D\,\sum_{s=r}^{\infty}\left((n+(s+1)(s+2))\,\hbox{\rm var}_s\psi + s\,\theta^{s} \right)
\right),\]
with $D:=\max(2,C)$. Since $\psi$ is H\"older continuous and $\theta \in (0,1)$, then
\[\epsilon_{r,n}:=D\,\sum_{s=r}^{\infty}
\left( (n+ (s+1)(s+2)) \hbox{\rm var}_s\psi+ s\,\theta^{s}\right)\rightarrow 0 \, 
\text{ when } \, r\to\infty,\]
for each $n,r\in {\mathbb N}$. We conclude that, 
$\mu[\boldsymbol{w}]:=\lim_{r\to\infty}\mu_{\psi_r}[\boldsymbol{w}]$ 
exists for each $\boldsymbol{w}\in\cup_{k=0}^\infty A^k$, and we have
\[\frac{\mu_{\psi_r}[\boldsymbol{w}]}{\mu[\boldsymbol{w}]}
\lessgtr
\exp\left(\pm\,D\,\sum_{s=r}^{\infty}\left((|\boldsymbol{w}|+(s+1)(s+2))\hbox{\rm var}_s\psi+ 
s\,\theta^{s}\right) \right),\]
for every $r\in{\mathbb N}$ and $\boldsymbol{w}\in \cup_{k=1}^\infty A^k$.

\medskip\noindent
\subsubsection{The limit $\lim_{r\to\infty}\mu_{\psi_r}$ is the Gibbs measure $\mu_\psi$.}\

\medskip\noindent It only remains to prove that $\mu$ such that 
$\mu[\boldsymbol{w}]:=\lim_{r\to\infty}\mu_{\psi_r}$ coincides with the original Gibbs measure 
$\mu_\psi$. Note first that $\mu$ so defined is $T$--invariant. Indeed, it is the weak${}^*$ 
limit of the sequence $(\mu_{\psi_r})_{r\geq 1}$ of $T$--invariant Markov approximants, it is a 
$T$--invariant probability measure as well.

\noindent Now, replacing $\psi_r$ by $\psi\pm\hbox{\rm var}_r\psi$, and making $p=(r+1)(r+2)+n-1$ 
in Ineq.~\eqref{periodic-approximation2}, it follows that
\begin{eqnarray*}
\mu[\boldsymbol{w}]  &\lessgtr &
\mu_{\psi_r}[\boldsymbol{w}]\,\exp\left(\pm\,\epsilon_{r,n}\right) \nonumber\\
                     &\lessgtr &
\frac{\sum_{\boldsymbol{a}\in\hbox{\rm \tiny Per}_p(A^{\mathbb N})\cap\,[\boldsymbol{w}]}
    e^{S_p\psi(\boldsymbol{a})}
     }{\sum_{\boldsymbol{a}\in\hbox{\rm \tiny Per}_p(A^{\mathbb N})}
     e^{S_p\psi(\boldsymbol{a})}}\,\exp\left(\pm\,2\epsilon_{r,n}\right)
\end{eqnarray*}
for every $\boldsymbol{w}\in \cup_{k=1}^nA^k$. By taking $n=r^2$, we obtain
\begin{equation}\label{periodic-approximation3}
\mu[\boldsymbol{w}]\lessgtr 
\frac{
\sum_{\boldsymbol{a}\in\hbox{\rm \tiny Per}_{(2r+1)(r+1)}(A^{\mathbb N})\cap\,[\boldsymbol{w}]}
e^{S_{(2r+1)(r+1)}\psi(\boldsymbol{a})}
     }{\sum_{\boldsymbol{a}\in\hbox{\rm \tiny Per}_{(2r+1)(r+1)}(A^{\mathbb N})}
e^{S_{(2r+1)(r+1)}\psi(\boldsymbol{a})}}\,\exp\left(\pm\,2\epsilon_{r,r^2}\right)
\end{equation}
for each $r\in{\mathbb N}$ and $\boldsymbol{w}\in \cup_{k=1}^{r^2}A^{k}$. On the other hand, 
the Gibbs measure $\mu_{\psi}$, whose existence is ensured by the fact that 
$\sum_{r}\hbox{\rm var}_r\psi<\infty$, is such that
\[\mu_\psi[\boldsymbol{w}]\lessgtr C^{\pm 1}\, 
\frac{\sum_{\boldsymbol{a}\in\hbox{\rm \tiny Per}_p(A^{\mathbb N})\cap\,[\boldsymbol{w}]}
e^{S_p\psi(\boldsymbol{a})} 
    }{\sum_{\boldsymbol{a}\in\hbox{\rm \tiny Per}_p(A^{\mathbb N})}e^{S_p\psi(\boldsymbol{a})}},\]
for each  $\boldsymbol{w}\in A^k$ with $k\leq p$. Since $\epsilon_{r,r^2}\rightarrow 0$ when
$r\to\infty$, it follows from this and Ineq.~\eqref{periodic-approximation3} that $\mu$ is 
absolutely continuous with respect to $\mu_\psi$. The Ergodic Decomposition Theorem implies
that $\mu_\psi$ is the only ergodic measure entering in the decomposition of the invariant 
measure $\mu$, therefore $\mu=\mu_{\psi}$.


\end{document}